\documentclass[11pt]{amsart}  

\usepackage[utf8]{inputenc}

\usepackage{graphicx,amssymb,amsmath,amsthm,wrapfig}
\usepackage{caption}
\usepackage{subcaption}
\usepackage{indentfirst}
\usepackage{cancel}
\usepackage{lipsum}
\usepackage{epstopdf}
\usepackage{epsfig}
\usepackage{comment}
\usepackage{enumerate,color}
\usepackage{empheq}
\usepackage{fullpage}
\usepackage{hyperref}
\usepackage{xcolor}
\usepackage{float}
\usepackage{enumitem}

\newcommand{\veps}{\varepsilon}

\everymath{\displaystyle}
\theoremstyle{plain}
\newtheorem{theorem}{Theorem}[section]
\newtheorem{corollary}[theorem]{Corollary}
\newtheorem{lemma}[theorem]{Lemma}
\newtheorem{definition}[theorem]{Definition}
\newtheorem{proposition}[theorem]{Proposition}
\newtheorem{remark}[theorem]{Remark}

\numberwithin{equation}{section}
\theoremstyle{definition}
\newtheorem{assumption}{Assumption}

\title{Emergence of Cuspons through Wave Breaking}
\author{Yunjoo Kim and Bongsuk Kwon}

\address{(Yunjoo Kim) Department of Mathematical Sciences, Ulsan National Institute of Science and Technology, Ulsan, 44919, Korea}
\email{gomuli3@unist.ac.kr}

\address{(Bongsuk Kwon) Department of Mathematical Sciences, Ulsan National Institute of Science and Technology, Ulsan, 44919, Korea}
\email{bkwon@unist.ac.kr}

\date{\today}

\begin{document}

\maketitle

\begin{abstract}
We study the emergence of cuspons through wave breaking, with the
Camassa--Holm equation as a principal model. For a broad open class of
smooth initial data, the first singularity occurs at a unique point and
has the sharp local H\"older regularity \(C^{3/5}\). We construct a
continuation of the Lagrangian flow through this singular time and
reconstruct from it a conservative Eulerian weak solution. We show that
the first singularity does not simply persist after breaking; instead,
it bifurcates immediately into an anticuspon--cuspon pair, each having
the sharp local H\"older regularity \(C^{2/3}\), with opposite one-sided
slope orientations. The two singular branches move apart from the
original breaking point, and the portion of the profile between them
reverses its local monotonicity. 
The proof is based on an extended inverse-slope variable that remains
regular through wave breaking. At the first singular time, this variable
has a nondegenerate quadratic zero, which splits immediately afterward
into two simple zeros. The corresponding degeneracy of the Lagrangian
flow map converts the quadratic first-breaking geometry into the
\(C^{3/5}\) profile and the post-breaking simple-zero geometry into the
\(C^{2/3}\) cuspon profiles. The generalized Hunter--Saxton family
provides an explicit realization of the same Lagrangian mechanism and
also allows us to track the subsequent interaction and annihilation of
the singular branches. Thus an anticuspon--cuspon pair emerges
dynamically from smooth initial data through wave breaking, rather than
being imposed through a singular traveling-wave ansatz.
\end{abstract}

\section{Introduction}
\label{sec:introduction}

Cusped singularities arise naturally in nonlinear wave equations, often
in the form of singular traveling waves or other special coherent
structures. A more fundamental dynamical question is whether such
profiles can be generated directly from smooth initial data. In
particular, one would like to understand whether wave breaking can
dynamically create cuspons, what local singularity is selected at the
first breaking time, and how the resulting cusped branches evolve
afterward. The purpose of this paper is to identify a mechanism by which
cuspons are generated through wave breaking and to describe the
associated transition in local singularity structure.

To investigate this mechanism, we first study the Camassa--Holm (CH)
equation on the real line,
%
\begin{equation}
\label{CH_main}
u_t-u_{xxt}+2\gamma u_x+3u u_x
=
2u_xu_{xx}+u u_{xxx},
\qquad
x\in\mathbb R,
\end{equation}
where \(u=u(x,t)\) is the unknown wave variable and \(\gamma\geq0\) is the
linear-dispersion parameter. Introduced by Camassa and Holm as an
integrable model for unidirectional shallow-water waves, the equation is
notable for its peaked solitary waves, or peakons \cite{CH}. It can also
be written in the nonlocal form
\begin{equation*}
    u_t+u u_x+p_x=0,
    \qquad
    p
    =
    (1-\partial_x^2)^{-1}
    \left(
        u^2+\frac12u_x^2+2\gamma u
    \right).
\end{equation*}
This formulation makes visible the competition between nonlinear
transport and nonlocal forcing that underlies the singularity formation
studied below.

A fundamental feature of \eqref{CH_main} is wave breaking: a smooth
solution may remain bounded while its slope becomes unbounded from below
in finite time. Blow-up criteria and the wave-breaking mechanism have
been studied extensively; see, among many works, \cite{BC,CE,LO}. The
breakdown of the classical solution does not, however, terminate the
dynamics. Global conservative and dissipative continuations have been
constructed by passing to Lagrangian variables and adjoining quantities
that remain regular when the Eulerian slope blows up
\cite{BCconservative,BCdissipative,HRperiodic}. In the case
\(\gamma=0\), Bressan--Chen--Zhang established uniqueness in a natural
class of conservative solutions \cite{BCZ}.

These continuation theories answer the fundamental question of whether
the Camassa--Holm flow can be continued beyond wave breaking. They do
not, however, determine the local Eulerian geometry generated by that
continuation. In particular, they do not identify the shape of the
profile immediately after the first breaking time, nor do they explain
whether the post-breaking evolution dynamically selects one of the
singular geometries known from special solutions. The present paper
addresses this complementary geometric question: what does the
conservative continuation actually look like near a generic first
breaking point?

The Camassa--Holm equation admits several singular coherent structures.
Lenells classified its weak traveling waves, including peakons, cuspons,
stumpons, and composite waves \cite{Lenells}. A peakon has a corner with
finite one-sided slopes, whereas a cuspon is continuous and has a
vertical tangent: its derivative diverges with opposite signs on the two
sides of the crest. Such solutions demonstrate that \eqref{CH_main}
admits cusped profiles. They arise, however, from singular traveling-wave
ansatzes or other special exact constructions and therefore do not
explain whether cuspons can be generated dynamically from smooth initial
data through wave breaking.

A complementary line of work concerns the generic structure of singular
conservative solutions. Li and Zhang identified stable local normal
forms for conservative Camassa--Holm dynamics, including generic
\(C^{2/3}\) singular curves and isolated \(C^{3/5}\) transition
singularities \cite{LiZhang}. Their theory describes the local
regularity patterns that occur in generic conservative solutions. In our
preceding work \cite{KKY}, we approached the problem from the formation
side: for a broad open class of smooth initial data, we proved that the
first wave-breaking event occurs at a unique point, is governed by a
stable self-similar regime, and has the sharp local H\"older regularity
\(C^{3/5}\).

The present paper determines what happens next. We continue a
nondegenerate first breaking point through the singular time and
reconstruct the resulting Eulerian profile. The \(C^{3/5}\) singularity
at the first breaking time does not persist as an isolated singularity.
Instead, it is a transition singularity from which two distinct singular
branches emerge immediately after breaking. These branches form an ordered
anticuspon--cuspon pair: the left branch is a cusped local minimum and
the right branch is a cusped local maximum. Each post-breaking
singularity has the sharp local H\"older regularity \(C^{2/3}\).
The two singular branches then separate and move away from the original
breaking point in opposite directions.

The bifurcation is encoded by the geometry of the extended Lagrangian
variables. At the first breaking time, the extended inverse-slope
variable has a nondegenerate quadratic zero at the first breaking label.
Immediately afterward, this quadratic zero splits into two simple zeros.
At the same time, the Lagrangian flow map degenerates quadratically at
the first breaking point and linearly at each post-breaking branch. This
change in degeneracy produces the transition
\[
    C^{3/5}
    \quad\longrightarrow\quad
    C^{2/3},
\]
and determines the opposite one-sided slope orientations of the cuspon
and anticuspon. Thus the post-breaking exponents and cusp geometry arise
from a single Lagrangian bifurcation mechanism.

This provides a geometric refinement of both classical wave-breaking
theory and conservative-solution theory. Classical results detect the
loss of \(C^1\)-regularity through slope blow-up while the amplitude
remains bounded. Conservative theories provide a continuation beyond
that loss of classical regularity. Our contribution is to identify the
precise local shape produced by the continuation. Near the first
breaking point, the post-breaking profile is not an unspecified weak
extension with unbounded slope: it has a canonical
anticuspon--cuspon structure, sharp H\"older exponents, prescribed
one-sided slope orientations, and a definite outward dynamics.

We construct this geometry by extending the inverse slope, the
Lagrangian Jacobian, and the Lagrangian profile through the singular
labels. The Eulerian solution is reconstructed only after these
quantities have been continued across the degenerate set. We then prove
that the reconstructed solution satisfies the derivative-energy balance
law and hence defines a conservative weak solution.
In the case \(\gamma=0\), it belongs to the class of conservative
solutions considered by Bressan--Chen--Zhang \cite{BCZ} and therefore
agrees with the unique solution in that class. 
Consequently, the anticuspon--cuspon pair is not
an ad hoc weak continuation or a singular profile imposed in advance;
it is generated dynamically from smooth initial data by the canonical
conservative evolution.

The generalized Hunter--Saxton family provides an explicit model of the
same mechanism. For every positive even parameter \(\lambda\), we derive
the corresponding first-breaking and post-breaking regularity exponents.
The explicit Lagrangian formulas allow the dynamics to be continued
globally and make it possible to track the subsequent interaction of the
cuspon and anticuspon branches. We identify a sharp criterion for their
finite-time annihilation. At the collision time, two
\(C^{\lambda/(\lambda+1)}\) branches merge into a single
\(C^{(2\lambda-1)/(2\lambda+1)}\) transition singularity. Thus the local
birth of an anticuspon--cuspon pair and its later annihilation are two
manifestations of the same degeneracy and unfolding of the Lagrangian
flow map.

\subsection{Main result for the Camassa--Holm equation}

We consider \eqref{CH_main} with initial datum 
\(  u(x,0)=u_0(x).\)
Let \(u\) denote the corresponding classical solution on its maximal
interval of classical existence, and let
\begin{equation*}
    \eta_t(a,t)=u(\eta(a,t),t),
    \qquad
    \eta(a,0)=a,
\end{equation*}
be the associated Lagrangian flow.
Let \(t_*\) denote the first breaking time.
Under the generic assumptions introduced below, we will show that
the first breaking occurs at an isolated Lagrangian label. %
To describe the local structure of the first breaking event, for each \(a\) sufficiently close to a first breaking label \(a_*\), let \(\tau(a)\) denote the first breaking time along the characteristic issued from \(a\). The construction and regularity of \(\tau\) are established in Lemma~\ref{CHlem:cusp-pair-no-theta}.
%
%

\begin{assumption}[Initial data and first-breaking structure]
\label{ass:CH-main}
The initial datum \(u_0\) is sufficiently smooth and spatially localized;
for convenience, we assume that
\(
    u_0\in H^5(\mathbb R).
\)
The corresponding classical solution has a finite first breaking time
\(t_*\), attained at a unique Lagrangian label \(a_*\). Moreover, the
local breaking-time map \(\tau\) is \(C^2\) near \(a_*\) and has a
nondegenerate minimum there:
\[
    \tau(a_*)=t_*,
    \qquad
    \tau'(a_*)=0,
    \qquad
    \tau''(a_*)>0.
\]
\end{assumption}

The Sobolev index \(5\) is not essential and is not expected to be
optimal. It is imposed only to ensure the classical regularity and the
Taylor expansions required in the analysis near the first breaking time.

Assumption~\ref{ass:CH-main} is satisfied by the admissible
class of smooth initial data constructed in \cite{KKY}. The precise
relation between this assumption and the sharp first-singularity structure
obtained in \cite{KKY} is explained in
Remark~\ref{rmk:CH-hypotheses}.

Under the above assumptions, we obtain the following description of the
first breaking event and its immediate post-breaking continuation.
%
%
%
%
 
\begin{theorem}[Post-breaking cusp formation for the Camassa--Holm equation]
\label{Thm_CH}
Suppose that Assumption~\ref{ass:CH-main} holds. Then the classical
solution $u$ reaches its first breaking time \(t_*\) at a unique Lagrangian
label \(a_*\). The classical Lagrangian flow $\eta$ admits a continuation to
\(
    0\leq t<t_*+\varepsilon_*
\)
for some \(\varepsilon_*>0\). The corresponding Eulerian reconstruction
defines a conservative weak solution on this interval and agrees with
the classical solution for \(0\leq t<t_*\).

At the first breaking time, the solution develops a single singularity
with exact local H\"older exponent \(3/5\). Immediately after \(t_*\),
this singularity bifurcates into an anticuspon--cuspon pair, each having
exact local H\"older exponent \(2/3\). More precisely, the following
assertions hold.

\begin{enumerate}
\item[\textup{(i)}]
The Lagrangian flow extends continuously to \(t=t_*\), and the solution
loses \(C^1\)-regularity at precisely one Eulerian point,
\(
    x_*:=\eta(a_*,t_*).
\) 
The profile has exact local H\"older exponent \(3/5\) at \(x_*\). Namely,
\(u(\cdot,t_*)\) belongs locally to \(C^{3/5}\) near \(x_*\), whereas
\(
    u(\cdot,t_*)\notin C^\alpha
\)
in any neighborhood of \(x_*\) for every \(\alpha>3/5\).

\item[\textup{(ii)}]
There exist two singular branches
\(
    x_-,x_+:
    (t_*,t_*+\varepsilon_*)\longrightarrow\mathbb R
\)
such that
\[
    x_-(t)<x_+(t),
    \qquad
    \lim_{t\downarrow t_*}x_-(t)
    =
    \lim_{t\downarrow t_*}x_+(t)
    =
    x_*.
\]
For every \(t\in(t_*,t_*+\varepsilon_*)\), the singular set emerging
locally from the first breaking point consists precisely of
\(x_-(t)\) and \(x_+(t)\).

At each of these points, the Eulerian profile has exact local H\"older
exponent \(2/3\). More precisely, \(u(\cdot,t)\) is locally \(C^{2/3}\)
near \(x_\pm(t)\), but does not belong locally to \(C^\alpha\) near either
point for any \(\alpha>2/3\).

The one-sided derivatives at \(x_+(t)\) satisfy
\begin{equation}
\label{CH_cuspon_slopes}
    \lim_{x\to x_+(t)^-}u_x(x,t)=+\infty,
    \qquad
    \lim_{x\to x_+(t)^+}u_x(x,t)=-\infty,
\end{equation}
whereas those at \(x_-(t)\) satisfy
\begin{equation}
\label{CH_anticuspon_slopes}
    \lim_{x\to x_-(t)^-}u_x(x,t)=-\infty,
    \qquad
    \lim_{x\to x_-(t)^+}u_x(x,t)=+\infty.
\end{equation}
Consequently, the profile has an anticuspon at \(x_-(t)\) and a cuspon
at \(x_+(t)\).
\end{enumerate}
\end{theorem}

\begin{figure}[t]
    \centering
    \includegraphics[width=0.93\textwidth]{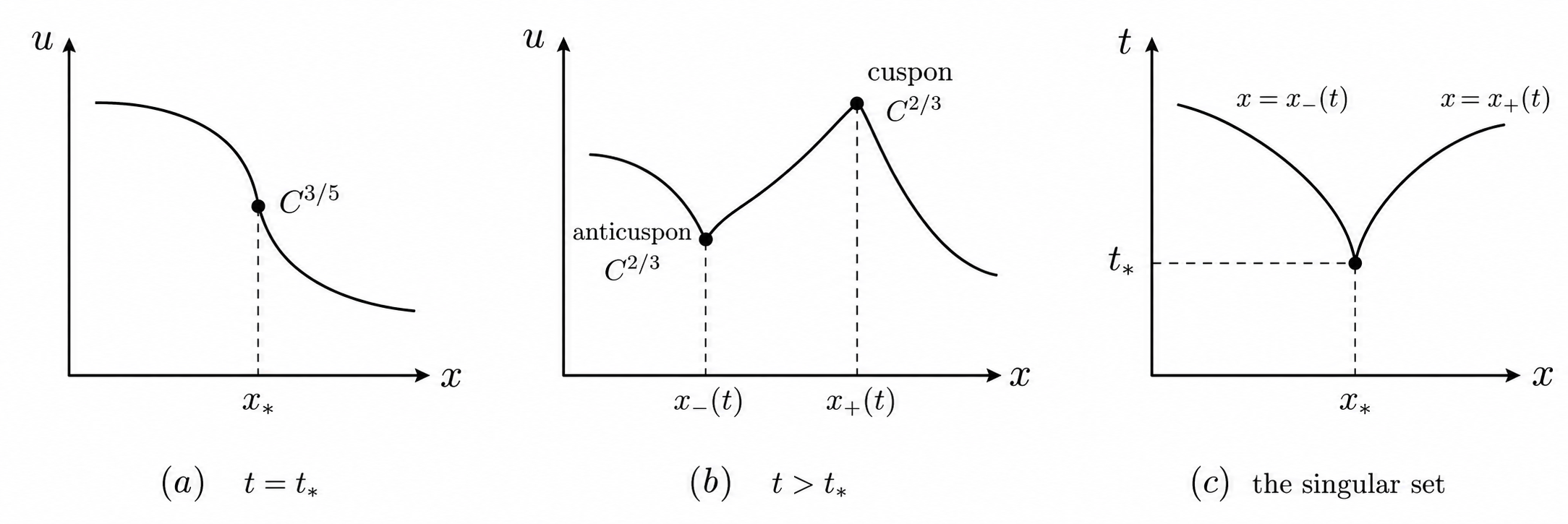}
    \caption{
    Schematic illustration of the post-breaking cusp formation described
    in Theorem~\ref{Thm_CH}.
    \textup{(a)} At the first breaking time \(t=t_*\), the Eulerian
    profile develops a single transition singularity at \(x_*\), with
    exact local H\"older exponent \(3/5\).
    \textup{(b)} For \(t>t_*\), this singularity bifurcates into an
    anticuspon at \(x_-(t)\) and a cuspon at \(x_+(t)\), each having
    exact local H\"older exponent \(2/3\).
    \textup{(c)} The corresponding singular set in the
    \((t,x)\)-plane. The two branches \(x=x_-(t)\) and \(x=x_+(t)\)
    emanate from the first breaking point \((t_*,x_*)\).
    The figure is schematic and is not drawn to scale.
    }
    \label{fig:CH-cuspon-bifurcation}
\end{figure}

The post-breaking geometry described in Theorem~\ref{Thm_CH} is
illustrated schematically in Figure~\ref{fig:CH-cuspon-bifurcation}.
The proof of Theorem~\ref{Thm_CH} is given in Section~\ref{sec3}.



\begin{remark}[On the hypotheses]
\label{rmk:CH-hypotheses}

\noindent
%
%
\textup{(i) Local nature of the nondegeneracy hypothesis.}
The essential local hypothesis for the post-breaking bifurcation is that
the first breaking label in
Assumption~\ref{ass:CH-main} is isolated and nondegenerate.
The same local argument may be applied separately near several isolated
nondegenerate first-breaking labels, provided that the corresponding
neighborhoods remain disjoint.

\medskip

\noindent
\textup{(ii) Relation with the first-singularity structure in
\cite{KKY}.}
The nondegeneracy condition on the breaking-time map \(\tau\) is the
Lagrangian counterpart of the sharp first-singularity structure obtained
in \cite{KKY}. 
We define the extended inverse-slope variable \(Q\) by
\(\textstyle 
    Q(a,t)=\frac{1}{u_x(\eta(a,t),t)}
\) 
away from wave breaking, and extend \(Q\) continuously through the
breaking set by setting \(Q=0\).
 As shown in
Corollary~\ref{cor:CH-quadratic-Q-expansion}, the extended inverse-slope
variable $Q$ 
 satisfies
\[
    Q(a,t_*)
    =
    -\frac14\tau''(a_*)(a-a_*)^2
    +
    o\bigl(|a-a_*|^2\bigr)
\]
as \(a\to a_*\).  In particular,
\[
    Q(a_*,t_*)=0,
    \qquad
    Q_a(a_*,t_*)=0,
    \qquad
    Q_{aa}(a_*,t_*)
    =
    -\frac12\tau''(a_*).
\]
Consequently,
\[
    \tau''(a_*)>0
    \quad\Longleftrightarrow\quad
    Q_{aa}(a_*,t_*)<0.
\]
Thus, a nondegenerate minimum of the breaking-time map is equivalent to
a nondegenerate quadratic maximum of the extended inverse-slope variable
at the first breaking label.


This condition is satisfied by the admissible class of smooth initial data
constructed in \cite{KKY}. Indeed, the self-similar analysis in \cite{KKY}
shows that, up to the first singular time, the spatial derivative \(u_x\)
has a unique minimum at the modulation center and that this minimum is
nondegenerate. More precisely, in the self-similar variables of \cite{KKY},
one has
\(
    U_y(0,s)=-2,
    U_{yy}(0,s)=0,
    U_{yyy}(0,s)>0,
\)
and the estimates away from \(y=0\) imply
\(
    U_y(y,s)>-2
    \text{ for } y\neq 0.
\)
Here \(U\) denotes the rescaled slope profile used in the self-similar
formulation of \cite{KKY}. 
Thus the first singularity occurs at a unique spatial point and is generated
by a strict nondegenerate minimum of \(u_x\).

By the Eulerian--Lagrangian correspondence established above, this unique
first-singularity point corresponds to a unique breaking label \(a_*\).
In terms of the extended inverse-slope variable, the nondegeneracy is
expressed as
\[
    Q(a_*,t_*)=0,
    \qquad
    Q_a(a_*,t_*)=0,
    \qquad
    Q_{aa}(a_*,t_*)<0.
\]
Using Corollary~\ref{cor:CH-quadratic-Q-expansion}, this is equivalent to
\[
    \tau'(a_*)=0,
    \qquad
    \tau''(a_*)>0.
\]
Therefore, Assumption~\ref{ass:CH-main} is satisfied by the
open class of smooth initial data constructed in \cite{KKY}, rather than
only by specially tuned examples.
This is precisely the Lagrangian counterpart of the sharp
\(C^{3/5}\) first-singularity structure obtained in \cite{KKY}.

\medskip

\noindent
\textup{(iii) Finite-time breaking.}
The requirement \(t_*<\infty\) in
Assumption~\ref{ass:CH-main} is guaranteed, for example, by
the condition
\[
    u_0'(a_0)<-\sqrt{2K},
    \qquad
    K
    :=
    \frac32\|u_0\|_{H^1}^2
    +
    3|\gamma|\|u_0\|_{H^1}.
\]
This condition is used only to guarantee finite-time wave breaking and
may be replaced by any other standard wave-breaking criterion. By
itself, however, it does not imply that the first breaking label is
unique or that the breaking-time map has a nondegenerate minimum. A proof
of this sufficient criterion is given in
Proposition~\ref{prop:CH-finite-breaking}.

\end{remark}

\begin{remark}[Identification with the canonical conservative solution]
\label{rmk_BCZ}
For the case of \(\gamma =0\) in \eqref{CH_main}, 
the Eulerian solution produced by the post-breaking Lagrangian construction
belongs to the conservative solution class specified in
Section~\ref{sec_u_unique} following \cite{BCZ}. The identification argument given there shows
that this continuation agrees with the canonical conservative
Camassa--Holm solution.

\end{remark}


Theorem~\ref{Thm_CH} gives a dynamical refinement of the first-singularity
result in \cite{KKY}. At \(t=t_*\), the profile has the sharp
pre-shock-type H\"older exponent \(3/5\). 
Although this singularity is cusp-like in its H\"older scaling, it does
not yet have the one-sided slope structure of either a cuspon or an
anticuspon described in
\eqref{CH_cuspon_slopes}--\eqref{CH_anticuspon_slopes}. 
Immediately after
the breaking time, the singular point bifurcates into a \(C^{2/3}\)
anticuspon and a \(C^{2/3}\) cuspon, respectively located at
\(x_-(t)\) and \(x_+(t)\). The one-sided limits identify not only the strength of these singularities
but also their distinct slope configurations, which become important in
the interaction analysis developed below. 

To isolate the underlying Lagrangian mechanism from the specific
nonlocal structure of the Camassa--Holm equation, we also study the
generalized Hunter--Saxton family, for which the same bifurcation
mechanism becomes explicit.

\subsection{The generalized Hunter--Saxton family}

The generalized Hunter--Saxton family is given by 
\begin{equation}
\label{genHS0}
    v_{tx}+v v_{xx}
    =
    -\frac1\lambda v_x^2,
    \qquad
    \lambda\in 2\mathbb{N}=\{2,4,6,\ldots\}.
\end{equation}
The standard Hunter--Saxton equation corresponds to \(\lambda=2\) and
was introduced as an asymptotic model for director-field dynamics
\cite{HS}. Generic regularity and Lipschitz metrics for conservative
solutions of Hunter--Saxton type equations were studied in
\cite{CaiChenShenTan}; see also \cite{GLW,GH} for further results on the
regularity and uniqueness of conservative solutions.

The restriction to even \(\lambda\) is essential for our continuation
mechanism. Indeed, the Jacobian of the extended Lagrangian flow is the \(\lambda\)-th power
of a quantity that changes sign at breaking. Consequently,
the extended flow remains monotone after its Jacobian vanishes, allowing
one to recover a single-valued Eulerian solution beyond the first
singular time.

For a classical solution \(v\) of \eqref{genHS0}, let
\[
    \zeta_t(a,t)=v(\zeta(a,t),t),
    \qquad
    \zeta(a,0)=a,
\]
denote the associated Lagrangian flow. It is convenient to introduce
the degeneracy factor
\[
    D_{\mathrm{HS}}(a,t)
    :=
    1+\frac{t}{\lambda}v_0'(a).
\]
As shown below, the singular labels are precisely the zeros of
\(D_{\mathrm{HS}}\).

\begin{assumption}[Initial data and first-breaking structure]
\label{ass:HS-main}
The initial datum satisfies the following conditions.

\begin{enumerate}
\item[\textup{(H1)}]
The initial profile is smooth and its derivative satisfies
\begin{equation*}
    v_0\in C^\infty(\mathbb R),
    \qquad
    v_0'\in L^1(\mathbb R)\cap L^\lambda(\mathbb R).
\end{equation*}

\item[\textup{(H2)}]
The initial slope \(v_0'\) attains a unique nondegenerate negative
minimum at a label \(a_0\):
\begin{equation*}
    v_0'(a_0)
    =
    \min_{a\in\mathbb R}v_0'(a)
    <0,
    \qquad
    v_0''(a_0)=0,
    \qquad
    v_0^{(3)}(a_0)>0.
\end{equation*}

\item[\textup{(H3)}]
The initial slope satisfies
\[
v_0'(a)\rightarrow 0
\qquad
\text{as }|a|\to\infty.
\]
\end{enumerate}
\end{assumption}
The decay condition \textup{(H3)} ensures that the far-field slope remains separated from the distinguished negative minimum \(v_0'(a_0)\), so that the first-breaking structure is localized near \(a_0\).

In both the Camassa--Holm and generalized Hunter--Saxton settings, we use
the notation \(a_\pm(t)\) for the two post-breaking singular labels and
\(x_\pm(t)\) for their corresponding Eulerian cusp locations. Since the
two equations are treated separately, no confusion should arise.

\begin{theorem}[Post-breaking cusp formation for the generalized
Hunter--Saxton equation]
\label{mainthm1}
Suppose that Assumption~\ref{ass:HS-main} holds. Then the first singular
time is
\(
    \textstyle t_*=-\frac{\lambda}{v_0'(a_0)},
\)
and the first breaking occurs at the unique Lagrangian label \(a_0\).

There exists \(\veps_{\mathrm{HS}}>0\) such that the classical
Lagrangian flow admits a continuation to
\(0\leq t<t_*+\veps_{\mathrm{HS}}\).
The corresponding Eulerian reconstruction defines a continuous function
\[
    v:
    \mathbb R\times[0,t_*+\veps_{\mathrm{HS}})
    \longrightarrow\mathbb R,
    \qquad
    v(\cdot,0)=v_0,
\]
such that, for every \(T<t_*+\veps_{\mathrm{HS}}\), \(v\) is an
Eulerian weak solution of \eqref{genHS0} on
\(\mathbb R\times(0,T)\).

Moreover, the exact conservation law
\begin{equation}
\label{mainthm_space'}
    \|v_x(\cdot,t)\|_{L^\lambda(\mathbb R)}
    =
    \|v_0'\|_{L^\lambda(\mathbb R)}
\end{equation}
holds for every \(0\leq t<t_*+\veps_{\mathrm{HS}}\). In particular,
\begin{equation}
\label{mainthm_space}
    v_x
    \in
    L^\infty
    \bigl(
        0,t_*+\veps_{\mathrm{HS}};
        L^\lambda(\mathbb R)
    \bigr),
    \qquad
    v
    \in
    L^\infty_{\mathrm{loc}}
    \bigl(
        \mathbb R\times[0,t_*+\veps_{\mathrm{HS}})
    \bigr).
\end{equation}

The local singularity structure is described as follows.

\begin{enumerate}
\item[\textup{(i)}]
At the first singular time \(t=t_*\), the solution loses
\(C^1\)-regularity at precisely one Eulerian point
\(x_*:=\zeta(a_0,t_*)\).
The pointwise H\"older exponent of \(v(\cdot,t_*)\) at \(x_*\) is
exactly
\(
    \textstyle \frac{2\lambda-1}{2\lambda+1}.
\)
More precisely,
\[
    0
    <
    \lim_{x\to x_*}
    \frac{
        |v(x,t_*)-v(x_*,t_*)|
    }{
        |x-x_*|^{\frac{2\lambda-1}{2\lambda+1}}
    }
    <
    \infty.
\]
In particular,
\(
    v(\cdot,t_*)
    \in
    C^{\frac{2\lambda-1}{2\lambda+1}}
\)
locally near \(x_*\), and this exponent is sharp.

\item[\textup{(ii)}]
For every \(t\in(t_*,t_*+\veps_{\mathrm{HS}})\), there exist exactly
two singular labels near \(a_0\),
\(
    a_-(t)<a_0<a_+(t),
\)
satisfying \(D_{\mathrm{HS}}(a_\pm(t),t)=0\). They bifurcate from
\(a_0\) in the sense that
\(
    a_\pm(t)\to a_0
\)
as \(t\downarrow t_*\).

Their Eulerian images
\(
    x_\pm(t):=\zeta(a_\pm(t),t)
\)
satisfy
\[
    x_-(t)<x_+(t),
    \qquad
    x_\pm(t)\longrightarrow x_*
    \quad\text{as }t\downarrow t_*.
\]
Thus the single first-breaking singularity bifurcates into a pair of
post-breaking singularities.

At each \(x_\pm(t)\), the pointwise H\"older exponent of
\(v(\cdot,t)\) is exactly
\(
    \textstyle \frac{\lambda}{\lambda+1}.
\)
More precisely,
\[
    0
    <
    \lim_{x\to x_\pm(t)}
    \frac{
        |v(x,t)-v(x_\pm(t),t)|
    }{
        |x-x_\pm(t)|^{\frac{\lambda}{\lambda+1}}
    }
    <
    \infty.
\]
In particular,
\(
    v(\cdot,t)\in C^{\frac{\lambda}{\lambda+1}}
\)
locally near \(x_\pm(t)\), and this exponent is sharp.

The one-sided slopes satisfy
\[
    \lim_{x\to x_+(t)^-}v_x(x,t)=+\infty,
    \qquad
    \lim_{x\to x_+(t)^+}v_x(x,t)=-\infty,
\]
whereas
\[
    \lim_{x\to x_-(t)^-}v_x(x,t)=-\infty,
    \qquad
    \lim_{x\to x_-(t)^+}v_x(x,t)=+\infty.
\]
Consequently, \(x_+(t)\) is a cuspon of local-maximum type, whereas
\(x_-(t)\) is an anticuspon of local-minimum type.

\item[\textup{(iii)}]
For every \(t\in(t_*,t_*+\veps_{\mathrm{HS}})\), the solution is
\(C^1\) away from the cusp pair:
\[
    v(\cdot,t)
    \in
    C^1
    \bigl(
        \mathbb R\setminus\{x_-(t),x_+(t)\}
    \bigr).
\]
Moreover, \(v_x\) is \(C^1\) in space-time away from the two singular
branches, and \(v\) is a classical solution of \eqref{genHS0} there.

\item[\textup{(iv)}]
The two cusp locations move outward from the central trajectory
\(x_0(t):=\zeta(a_0,t)\). More precisely,
\[
    x_-'(t)<x_0'(t)<x_+'(t).
\]
\end{enumerate}
\end{theorem}

We remark that 
\(\textstyle 
    \frac{\lambda}{\lambda+1}
    >
    \frac{2\lambda-1}{2\lambda+1},
\)
so the two post-breaking cuspons have strictly higher pointwise
H\"older regularity than the profile at the first breaking time.
When \(\lambda=2\), the sharp exponents are \(3/5\) at the first
breaking time and \(2/3\) at each post-breaking cuspon.

\begin{remark}[On the uniqueness of the first breaking label]
\label{rmk:HS-uniqueness}
The uniqueness assumption in
Assumption~\ref{ass:HS-main}\textup{(H2)}
is imposed only to exclude simultaneous first-breaking events.
If \(v_0'\) has finitely many isolated nondegenerate global minima
with the same negative value, then the local post-breaking construction
and cusp-bifurcation analysis can be carried out independently near each
corresponding first-breaking label, provided that the associated
neighborhoods remain disjoint.
\end{remark}

\subsection{Global anticuspon--cuspon dynamics}
	Starting from the time extension of the pre-breaking Lagrangian flow
	obtained in Theorem~\ref{mainthm1}, the simple form of \eqref{genHS0} enables us
	to derive an explicit global-in-time representation of \(v\).
This representation shows that the cuspon-pair formation described in
Theorem~\ref{mainthm1} is not confined to the first breaking event:
beyond the local time interval \([0,t_*+\veps_{\mathrm{HS}})\),
further anticuspon--cuspon pairs may be generated at other
nondegenerate negative local minima of \(v_0'\).
Conversely, a negative local maximum of \(v_0'\) may lead to the reverse
phenomenon: an adjacent inner cuspon and anticuspon can collide and
annihilate.


We impose the following global condition on the critical points of the
initial slope.

\begin{assumption}[Global critical-point structure]
\label{ass:HS-global}
Define
\[
    \mathcal A
    :=
    \{a\in\mathbb R:\ v_0''(a)=0\},
    \qquad
    \mathcal A^-
    :=
    \{a\in\mathcal A:\ v_0'(a)<0\}.
\]
We assume the following.

\begin{enumerate}
\item[\textup{(G1)}]
The set \(\mathcal A^-\) is locally finite.

\item[\textup{(G2)}]
Every negative critical point of \(v_0'\) is nondegenerate:
\[
    v_0^{(3)}(a)\neq0
    \qquad
    \text{for every }a\in\mathcal A^-.
\]
\end{enumerate}
\end{assumption}

\begin{theorem}[Global anticuspon--cuspon dynamics for the generalized
Hunter--Saxton equation]
\label{thm_collapse}
Suppose that Assumptions~\ref{ass:HS-main} and
\ref{ass:HS-global} hold. Then the weak solution constructed in
Theorem~\ref{mainthm1} extends globally in time. Moreover, the conservation
law
\[
    \|v_x(\cdot,t)\|_{L^\lambda(\mathbb R)}
    =
    \|v_0'\|_{L^\lambda(\mathbb R)}
\]
holds for every \(t\geq0\).

Let \(a'<a''\) be two consecutive nondegenerate negative local minima of
\(v_0'\), and assume that \(v_0'\) has a unique critical point
\(\bar a\in(a',a'')\). Thus \(\bar a\) is a nondegenerate local maximum,
and
\[
    v_0^{(3)}(a')>0,
    \qquad
    v_0^{(3)}(\bar a)<0,
    \qquad
    v_0^{(3)}(a'')>0.
\]
Consider the right-hand cuspon branch generated at \(a'\) and the
left-hand anticuspon branch generated at \(a''\). Then the following
alternative holds.

\begin{enumerate}
\item[\textup{(i)}]
If
\(
    v_0'(\bar a)<0,
\)
then the two inner branches meet at the finite annihilation time
\(\textstyle 
    t_{\mathrm{ann}}
    :=
    -\frac{\lambda}{v_0'(\bar a)}.
\)
At \(t=t_{\mathrm{ann}}\), their Eulerian locations and profile values
coincide, and the two singular branches coalesce into a single singular
point. For times immediately after \(t_{\mathrm{ann}}\), the two inner
branches are absent; in this sense, the anticuspon--cuspon pair
annihilates.

Each incoming branch has exact local H\"older exponent
\( \textstyle 
    \frac{\lambda}{\lambda+1},
\)
whereas the profile at the collision point has exact local H\"older
exponent
\( \textstyle
    \frac{2\lambda-1}{2\lambda+1}.
\)
Thus the annihilation event produces an instantaneous loss of Eulerian
H\"older regularity.

\item[\textup{(ii)}]
If
\(
    v_0'(\bar a)\geq0,
\)
then the two inner branches do not meet at any finite positive time.
In particular, no finite-time anticuspon--cuspon annihilation occurs
between the branches generated at \(a'\) and \(a''\).
\end{enumerate}
\end{theorem}

Thus, the sign of \(v_0'(\bar a)\) determines whether the two inner
branches undergo annihilation in finite positive time. 


\subsection{Organization of the paper}
Section~\ref{sec3} is devoted to the Camassa--Holm equation.  In
Sections~\ref{sec_CH_construct}--\ref{sec_CH_reg}, we construct the
post-breaking Lagrangian continuation, prove the emergence of the
anticuspon--cuspon pair, and determine the sharp local regularity. 
Section~\ref{sec_u_unique} verifies the Eulerian weak formulation and
identifies the constructed continuation with the unique conservative
solution in the class of \cite{BCZ}, completing the proof of
Theorem~\ref{Thm_CH}.

Section~\ref{sec2} treats the generalized Hunter--Saxton equation.  We
construct the post-breaking Eulerian weak solution, prove the conservation
law \eqref{mainthm_space'}, and establish the local cusp formation and sharp
regularity asserted in Theorem~\ref{mainthm1}.  In
Section~\ref{sec_HS_global}, we extend this construction globally, analyze
the interaction and annihilation of singular branches, and prove
Theorem~\ref{thm_collapse}.

\section{The Camassa--Holm equation}\label{sec3}
In this section, we prove Theorem~\ref{Thm_CH}.
The proof is developed through a series of lemmas and propositions. 
We now return to the Camassa--Holm equation \eqref{CH_main} on
\(\mathbb R\), supplemented with the initial condition
\begin{equation}\label{CH_init}
    u(x,0)=u_0(x).
\end{equation}
Throughout this section, \(u\) denotes the classical solution on its
maximal interval of classical existence.  Equation \eqref{CH_main} is
 equivalently written as the nonlocal system
\begin{subequations}\label{CH_nonlocal}
\begin{align}
    u_t+u u_x&=-p_x,
    \label{CH_transport}
    \\
    p-p_{xx}&=u^2+\frac12u_x^2+2\gamma u,
    \label{CH_elliptic}
\end{align}
\end{subequations}
where \(p=p(x,t)\) denotes the Eulerian nonlocal potential.  With
\[
    G(x):=\frac12e^{-|x|},
    \qquad
    (1-\partial_x^2)G=\delta_0,
\]
we have
\begin{equation}\label{CH_potential_representation}
    p
    =
    G*\left(u^2+\frac12u_x^2+2\gamma u\right).
\end{equation}

\subsection{Energy conservation and bounds for the nonlocal potential}

\begin{lemma}[Conservation of the \(H^1\)-energy]
For every smooth solution of \eqref{CH_main} for which the integrations
below are justified,
\begin{equation}\label{CH_cons_smooth}
    \int_{\mathbb R}\bigl(u^2+u_x^2\bigr)(x,t)\,dx
    =
    \int_{\mathbb R}\bigl(u_0^2+(u_0')^2\bigr)(x)\,dx
\end{equation}
throughout the classical lifespan.
\end{lemma}

\begin{proof}
Multiplying \eqref{CH_main} by \(u\), integrating over \(\mathbb R\),
and integrating by parts, we obtain
\[
    \int_{\mathbb R}u u_t\,dx
    -\int_{\mathbb R}u u_{xxt}\,dx
    =
    \frac12\frac{d}{dt}
    \int_{\mathbb R}\bigl(u^2+u_x^2\bigr)\,dx.
\]
The remaining terms vanish:
\[
    \int_{\mathbb R}2\gamma u u_x\,dx=0,
    \qquad
    \int_{\mathbb R}3u^2u_x\,dx=0,
\]
and
\[
    \int_{\mathbb R}
    \left(2u u_xu_{xx}+u^2u_{xxx}\right)\,dx
    =
    \int_{\mathbb R}(u^2u_{xx})_x\,dx
    =0.
\]
This proves \eqref{CH_cons_smooth}.
\end{proof}

\begin{lemma}[Bounds for the nonlocal potential]
\label{CH_lem_up}
Let \(u\in H^1(\mathbb R)\), and let \(p\) be determined by
\eqref{CH_elliptic}.  Then $\|u\|_{L^\infty}\le \|u\|_{H^1}$, \(p\in H^1(\mathbb R)\), and
\begin{equation}\label{CH_p_Linfty_bound}
    \max\bigl\{
        \|p\|_{L^\infty},
        \|p_x\|_{L^\infty}
    \bigr\}
    \leq
    \frac12\|u\|_{H^1}^2
    +|\gamma|\|u\|_{L^2}.
\end{equation}
Moreover,
\begin{equation}\label{CH_p_L2_bound}
    \max\bigl\{
        \|p\|_{L^2},
        \|p_x\|_{L^2}
    \bigr\}
    \leq
    \frac12\|u\|_{H^1}^2
    +2|\gamma|\|u\|_{L^2}.
\end{equation}
\end{lemma}

\begin{proof}
By \eqref{CH_potential_representation},
\[
    p
    =
    G*\left(u^2+\frac12u_x^2+2\gamma u\right),
    \qquad
    p_x
    =
    G_x*\left(u^2+\frac12u_x^2+2\gamma u\right).
\]
Recall
\[
    \|G\|_{L^\infty}
    =\|G_x\|_{L^\infty}
    =\frac12,
    \qquad
    \|G\|_{L^2}
    =\|G_x\|_{L^2}
    =\frac12,
\]
\[
    \|G\|_{L^1}
    =\|G_x\|_{L^1}
    =1,
    \qquad
    \left\|u^2+\frac12u_x^2\right\|_{L^1}
    \leq
    \|u\|_{H^1}^2.
\]
For $K\in\{G,G_x\}$, Young's convolution inequality gives, for the $L^\infty$ bound,
\begin{equation*}
	\left\|K\ast \left(u^2+\frac12u_x^2+2\gamma u\right)\right\|_{L^{\infty}}\leq \|K\|_{L^{\infty}}\left\|u^2+\frac12u_x^2\right\|_{L^1}+\|K\|_{L^2}\|2\gamma u\|_{L^2}\leq \frac{1}{2}\|u\|_{H^1}^2+|\gamma|\|u\|_{L^2},
\end{equation*}
which yields \eqref{CH_p_Linfty_bound}. For the $L^2$ bound, 
\begin{equation*}
	\left\|K\ast \left(u^2+\frac12u_x^2+2\gamma u\right)\right\|_{L^2}\leq \|K\|_{L^2}\left\|u^2+\frac12u_x^2\right\|_{L^1}+\|K\|_{L^1}\|2\gamma u\|_{L^2}\leq \frac{1}{2}\|u\|_{H^1}^2+2|\gamma|\|u\|_{L^2},
\end{equation*}
which yields \eqref{CH_p_L2_bound}.
\end{proof}

\subsection{Lagrangian slope dynamics and the inverse-slope variable}\label{sec:CH-inverse-slope}

Let \(\eta\) denote the Lagrangian flow of the classical solution:
\begin{equation}\label{CH_flow}
    \eta_t(a,t)=u(\eta(a,t),t),
    \qquad
    \eta(a,0)=a.
\end{equation}
We introduce the Lagrangian variables
\begin{equation*}
\begin{aligned}
    U(a,t)&:=u(\eta(a,t),t),\\
    P(a,t)&:=p(\eta(a,t),t),\\
    S(a,t)&:=u_x(\eta(a,t),t),\\
    J(a,t)&:=\eta_a(a,t).
\end{aligned}
\end{equation*}

Differentiating \eqref{CH_transport} with respect to \(x\), using
\eqref{CH_elliptic}, and then restricting to the flow \eqref{CH_flow}, we
obtain
\begin{equation}\label{CH_slope_characteristic}
    \partial_t S
    =
    -\frac12S^2
    +\mathcal F,
\end{equation}
where
\begin{equation}\label{F_def}
    \mathcal F(a,t)
    :=
    U(a,t)^2+2\gamma U(a,t)-P(a,t).
\end{equation}
Also, by \eqref{CH_flow}, the Lagrangian Jacobian satisfies
\begin{equation*}
    \partial_tJ
    =
    SJ,
    \qquad
    J(a,0)=1.
\end{equation*}

By Lemma~\ref{CH_lem_up} and the conservation law
\eqref{CH_cons_smooth}, the forcing in
\eqref{CH_slope_characteristic} is uniformly bounded throughout the
classical lifespan:
\begin{equation}\label{CH_forcing_bound}
    |\mathcal F(a,t)|
    \leq K,
    \qquad
    K
    :=
    \frac32\|u_0\|_{H^1}^2
    +3|\gamma|\|u_0\|_{H^1}.
\end{equation}

\begin{proposition}[A sufficient condition for finite-time wave breaking]
\label{prop:CH-finite-breaking}
Let \(a_0\in\mathbb R\) satisfy
\begin{equation*}
    u_0'(a_0)<-\sqrt{2K},
\end{equation*}
where \(K\) is defined in \eqref{CH_forcing_bound}.  Then the classical
solution undergoes wave breaking in finite time.  More precisely,
\(S(a_0,t)\) becomes unbounded from below no later than
\begin{equation*}
    T_{\mathrm{cmp}}
    :=
    \frac{1}{\sqrt{2K}}
    \log\left(
        \frac{u_0'(a_0)-\sqrt{2K}}
             {u_0'(a_0)+\sqrt{2K}}
    \right).
\end{equation*}
\end{proposition}

\begin{proof}
Let \(g\) solve
\begin{equation}\label{eq_g}
    g'(t)=-\frac12g(t)^2+K,
    \qquad
    g(0)=u_0'(a_0).
\end{equation}
Since \(\mathcal F\leq K\), comparison with
\eqref{CH_slope_characteristic} gives
\(
    S(a_0,t)\leq g(t)
\)
for as long as both quantities remain finite.  Set
\( \textstyle 
    \kappa:=\sqrt{2K},
    C_0:=\frac{g(0)-\kappa}{g(0)+\kappa}>1.
\)
The solution of \eqref{eq_g} is
\[
    g(t)
    =
    \kappa
    \frac{1+C_0e^{-\kappa t}}
         {1-C_0e^{-\kappa t}}.
\]
Consequently,
\[
    g(t)\longrightarrow-\infty
    \qquad
    \text{as}
    \qquad
    t\nearrow
    \frac{1}{\kappa}\log C_0
    =T_{\mathrm{cmp}}.
\]
It follows that \(S(a_0,t)\), and hence \(u_x\), becomes
unbounded from below in finite time.
\end{proof}

\begin{remark}[Wave-breaking criteria for the Camassa--Holm equation]
\label{rmk:CH-wave-breaking-literature}
The condition in Proposition~\ref{prop:CH-finite-breaking} is not meant
to be optimal; it provides a simple sufficient criterion for the
post-breaking analysis developed below. The wave-breaking theory for
the Camassa--Holm equation is considerably more extensive. Various
sufficient conditions, including local-in-space criteria, are known,
as well as a necessary and sufficient characterization of finite-time
wave breaking. We refer to
\cite{BrandoleseLocal,BC, CE,McKean}
and the references therein.
\end{remark}

\subsection{Post-breaking Lagrangian continuation}
\label{sec_CH_construct}

In this section, we construct solution in following steps: a regular
coupled system, its local continuation, positivity and monotonicity, a
regularity bootstrap, and the reconstruction of the Eulerian profile.

\subsubsection{The extended coupled Lagrangian system}
\label{sec:CH-extended-coupled-system}
The Eulerian slope \(S\) becomes unbounded at wave breaking, while the
classical Jacobian \(J\) may vanish.  We therefore rewrite the dynamics
in terms of Lagrangian variables that remain finite at the first
singular time.  For this purpose, we define
\[
    N:=SJ,
    \qquad
    M:=S^2J.
\]
We claim that the variables \(J,N,M\) remain finite even when \(S\) becomes unbounded.

It is convenient to work first with the cumulative energy variable
\begin{equation*}
    H(a,t)
    :=
    \int_{-\infty}^{a}
    \bigl(U(b,t)^2J(b,t)+M(b,t)\bigr)\,db.
\end{equation*}
At \(t=0\),
\begin{equation*}
    H(a,0)
    =
    \int_{-\infty}^{a} 
        | u_0|^2+ | u_0' |^2
    \,dx.
\end{equation*}

Using the Green function \(G(x)=\frac12e^{-|x|}\), define
\begin{align}
    P(a,t)
    &:=
    \int_{\mathbb R}
    G\bigl(\eta(a,t)-\eta(b,t)\bigr)
    \left[
        \frac12H_a(b,t)
        +
        \left(
            \frac12U(b,t)^2+2\gamma U(b,t)
        \right)
        \eta_a(b,t)
    \right]\,db,
    \label{CH_extended_P}
    \\
    \Pi(a,t)
    &:=
    \int_{\mathbb R}
    G_x\bigl(\eta(a,t)-\eta(b,t)\bigr)
    \left[
        \frac12H_a(b,t)
        +
        \left(
            \frac12U(b,t)^2+2\gamma U(b,t)
        \right)
        \eta_a(b,t)
    \right]\,db.
    \label{CH_extended_Pi}
\end{align}
On the classical lifespan, one can see by \eqref{CH_elliptic} that the change of variables
\(x=\eta(b,t)\) gives
\[
    P(a,t)=p(\eta(a,t),t),
    \qquad
    \Pi(a,t)=p_x(\eta(a,t),t).
\]

Recalling that
\[
    H_a=U^2J+M,
    \qquad
    H_t=U^3+2\gamma U^2-2UP,
\]
and using the commutation of the \(a\)- and \(t\)-derivatives on the
classical lifespan, we obtain
\[
    \partial_t\bigl(U^2J+M\bigr)
    =
    \partial_t H_a
    =
    \partial_a H_t
    =
    \partial_a
    \left(
        U^3+2\gamma U^2-2UP
    \right).
\]

Consequently, the regular variables \((\eta,U,H)\) satisfy
\begin{subequations}\label{CH_regular_lagrangian_system}
\begin{align}
    \eta_t&=U,
    \label{CH_regular_eta}
    \\
    U_t&=-\Pi,
    \label{CH_regular_U}
    \\
    H_t&=U^3+2\gamma U^2-2UP.
    \label{CH_regular_H}
\end{align}
\end{subequations}
Unlike the equations written in terms of the Eulerian slope
\[
    S=\frac{U_a}{\eta_a},
\]
the system \eqref{CH_regular_lagrangian_system} contains no division by
the Lagrangian Jacobian \(J=\eta_a\). Indeed, the nonlocal quantities
\(P\) and \(\Pi\) are defined directly from \(\eta\), \(U\), and \(H\). 
These expressions remain well-defined when \(J=\eta_a\) vanishes.
At such a label, the flow map may lose local invertibility, but the
variables \(\eta\), \(U\), and \(H\), as well as the Lagrangian energy
density
\(
    H_a=U^2J+M,
\)
remain finite. In particular, when \(J=0\), one has
\(
    H_a=M,
\)
so the derivative energy is still recorded by \(H\), even though the
Eulerian slope \(S=N/J\) may become infinite. Thus
\eqref{CH_regular_lagrangian_system} has no algebraic singularity at a
breaking label. The local well-posedness and preservation of the
admissible Lagrangian class are established in the next proposition.

Differentiating \eqref{CH_regular_lagrangian_system} with respect to
\(a\), 
one obtains
\begin{subequations}\label{CH_derived_JNM_system}
\begin{align}
    J_t&=N,
    \label{CH_derived_J}
    \\
    N_t&=\frac12M+\mathcal FJ,
    \label{CH_derived_N}
    \\
    M_t&=2\mathcal FN.
    \label{CH_derived_M}
\end{align}
\end{subequations}
Moreover,
\begin{equation}\label{CH_compatibility_H}
    N^2=JM
\end{equation}
is preserved by the evolution, because
\(
    \partial_t(N^2-JM)=0.
\)
The initial data are
\begin{equation}\label{CH_extended_initial_derivatives}
    J(a,0)=1,
    \qquad
    N(a,0)=u_0'(a),
    \qquad
    M(a,0)=\bigl(u_0'(a)\bigr)^2.
\end{equation}

We now introduce the inverse-slope variables in a local compressive
interval \(I\) containing the first breaking label and satisfying
\(u_0'(a)\neq0\) on \(I\).  Once \(\mathcal F\) is known from \eqref{F_def} with the variables solving 
\eqref{CH_regular_lagrangian_system}, define \(Q\) and \(R\) on \(I\)
as the solutions of
\begin{equation}\label{CH_Q_evolution}
    Q_t
    =
    \frac12-\mathcal FQ^2,
    \qquad
    Q(a,0)=\frac{1}{u_0'(a)},
\end{equation}
and
\begin{equation}\label{CH_R_Q_evolution}
    R_t
    =
    \bigl(1+2\mathcal F\bigr)
    \frac{Q}{1+Q^2}R,
    \qquad
    R(a,0)=1+\bigl(u_0'(a)\bigr)^2.
\end{equation}
Define
\begin{equation}\label{CH_JNM_QR_formulas}
    \widehat J
    :=
    \frac{Q^2R}{1+Q^2},
    \qquad
    \widehat N
    :=
    \frac{QR}{1+Q^2},
    \qquad
    \widehat M
    :=
    \frac{R}{1+Q^2}.
\end{equation}
A direct calculation shows that
\((\widehat J,\widehat N,\widehat M)\) satisfies
\eqref{CH_derived_JNM_system}, the compatibility relation
\(\widehat N^2=\widehat J\widehat M\), and the initial data
\eqref{CH_extended_initial_derivatives}.  By uniqueness for the
linearized system, it follows that
\[
    J=\widehat J,
    \qquad
    N=\widehat N,
    \qquad
    M=\widehat M
\]
on \(I\).  In particular,
\begin{equation}\label{CH_J_Q_R}
    J
    =
    \frac{Q^2R}{1+Q^2}, \qquad U_a
    =
    N
    =
    \frac{QR}{1+Q^2}, \qquad M
    =
    \frac{R}{1+Q^2}.
\end{equation}
Where \(J>0\), these formulas agree with
\[
    Q=\frac1S,
    \qquad
    R=(1+S^2)J.
\]
Thus \eqref{CH_Q_evolution}--\eqref{CH_R_Q_evolution} provide a regular
local description of the slope dynamics through \(Q=0\).

\subsubsection{Local continuous construction}
\label{sec:CH-local-continuous-construction}

\begin{proposition}[Local continuous continuation]
\label{prop:CH-local-continuation}
Suppose that Assumption~\ref{ass:CH-main} holds. Then there exist
\(\varepsilon_*>0\) and an extended Lagrangian solution
\(
    (\eta,U,H)
\)
of \eqref{CH_regular_lagrangian_system} on
\(
    0\leq t<t_*+\varepsilon_*.
\)
This solution agrees with the classical Lagrangian solution on
\(0\leq t<t_*\).  
Moreover, \(\eta\), \(U\), \(H\), \(P\), and \(\Pi\) are continuous
in \((a,t)\), the map \(\eta(\cdot,t)\) is nondecreasing
for every \(t\), and
\[
    \eta(\cdot,t)-\operatorname{Id}\in L^\infty(\mathbb R)
\]
on every bounded time interval.

The differentiated variables \(J,N,M\) satisfy
\eqref{CH_derived_JNM_system} and the compatibility relation
\eqref{CH_compatibility_H}.

Let \(I\) be a sufficiently small open interval containing the first
breaking label \(a_*\), chosen so that the classical flow is compressive
on \(I\) before the first breaking time. Then the variables \(Q\) and
\(R\), determined by
\eqref{CH_Q_evolution}--\eqref{CH_R_Q_evolution}, admit continuous
extensions to
\(
    I\times[0,t_*+\varepsilon_*)
\) 
and satisfy \eqref{CH_J_Q_R}.
\end{proposition}

%

\begin{proof}
Set
\(
    \widetilde{\eta}(a,t):=\eta(a,t)-a,
\)
and introduce
\[
    V
    :=
    \left\{
        f\in C_b(\mathbb R):
        f_a\in L^2(\mathbb R)
    \right\},
    \qquad
    \|f\|_V
    :=
    \|f\|_{L^\infty}+\|f_a\|_{L^2}.
\]
We regard
\(
    X(t)
    :=
    \bigl(
        \widetilde{\eta}(\cdot,t),
        U(\cdot,t),
        H(\cdot,t)
    \bigr)
\)
as an element of the Lagrangian state space
\(
    \mathcal E
    :=
    V\times H^1(\mathbb R)\times V.
\) 
Here and below,
\(
    \eta=\operatorname{Id}+\widetilde{\eta}
\)
in the definitions of \(P\) and \(\Pi\).

On bounded subsets of the admissible Lagrangian class, the maps
\[
    (\widetilde{\eta},U,H)
    \longmapsto
    (P,\Pi),
\]
defined by \eqref{CH_extended_P}--\eqref{CH_extended_Pi} with
\(\eta=\operatorname{Id}+\widetilde{\eta}\), are locally Lipschitz
with respect to the topology of \(\mathcal E\). Equivalently, on a
neighborhood of each bounded subset of the admissible class, the
regular Lagrangian system can be written as the Banach-space ordinary
differential equation
\[
    X_t=\mathcal F_{\mathrm{Lag}}(X),
\]
where
\[
    \mathcal F_{\mathrm{Lag}}
    (\widetilde{\eta},U,H)
    =
    \left(
        U,\,
        -\Pi,\,
        U^3+2\gamma U^2-2UP
    \right).
\]
Hence the Picard--Lindel\"of theorem yields a unique local solution.

	To continue this solution through the first breaking time $t_*$, we use the conservation of the total Lagrangian energy:
\[
    H(+\infty,t)
    =
    H(+\infty,0)
    =
    \|u_0\|_{H^1}^2.
\]
Consequently, the standard a priori estimates for the regular
Lagrangian system yield uniform bounds for the state variable
\(
    \bigl(
        \widetilde{\eta}(\cdot,t),
        U(\cdot,t),
        H(\cdot,t)
    \bigr)
\)
in \(\mathcal E\) on every bounded time interval on which the solution
exists. The corresponding estimates for the nonlocal quantities
\(P\) and \(\Pi\) remain controlled by the conserved energy. In
particular, the vector field
\(\mathcal F_{\mathrm{Lag}}\) remains bounded and locally Lipschitz
along the solution.

At the Eulerian breaking time \(t_*\), the Jacobian
\(
    J=\eta_a
\)
may vanish. However, the regular Lagrangian system
\eqref{CH_regular_lagrangian_system} contains no division by \(J\).
Thus the vanishing of \(J\) produces no algebraic singularity in the
regular Lagrangian variables. Moreover, the nonlocal quantities
\(P\) and \(\Pi\) remain controlled by the conserved energy.
Therefore the standard continuation criterion for the Banach-space
ordinary differential equation yields an extension of the regular
Lagrangian solution to
\(
    0\leq t<t_*+\varepsilon_*
\) 
for some \(\varepsilon_*>0\).

For \(0\leq t<t_*\), the regular Lagrangian system is obtained from
the classical Eulerian solution, and
\(\eta(\cdot,t)\) is a diffeomorphism. Uniqueness of the
Banach-space ordinary differential equation therefore implies that
the extended variables agree with the classical Lagrangian variables
on
\(
    0\leq t<t_*.
\)

Differentiating \eqref{CH_regular_lagrangian_system} with respect to
\(a\) yields \eqref{CH_derived_JNM_system}. Moreover, along this
system, 
\[
    \partial_t\bigl(N^2-JM\bigr)=0.
\]
Since the initial data satisfy $N(a,0)^2=J(a,0)M(a,0)$, it follows that $N^2=JM$ is preserved throughout the extended lifespan. Hence the
compatibility relation \eqref{CH_compatibility_H} holds on
\(
    0\leq t<t_*+\varepsilon_*.
\)

Finally, once the continuous forcing \(\mathcal F\) is known, for
each fixed \(a\in I\) the parameter-dependent ordinary differential
equations
\eqref{CH_Q_evolution}--\eqref{CH_R_Q_evolution}
have unique continuous solutions on
\(
    0\leq t<t_*+\varepsilon_*.
\)
The algebraic variables
\(
    (\widehat J,\widehat N,\widehat M)
\)
defined by \eqref{CH_JNM_QR_formulas} satisfy the same system
\eqref{CH_derived_JNM_system} and the same initial data as
\((J,N,M)\). By uniqueness,
\(
    (\widehat J,\widehat N,\widehat M)
    =
    (J,N,M),
\)
which proves \eqref{CH_J_Q_R}.
\end{proof}


Proposition~\ref{prop:CH-local-continuation} gives a continuous
post-breaking continuation.  The higher spatial regularity needed for
the cusp analysis will be obtained below by a separate bootstrap
argument.

\subsubsection{Positivity and monotonicity}
\label{sec:CH-positivity-monotonicity}

Let \(I_*\Subset I\) be a bounded interval containing \(a_*\), fix
\(0<\varepsilon_0<\varepsilon_*\), and set
\[
    \mathcal K
    :=
    \overline{I_*}
    \times
    [0,t_*+\varepsilon_0].
\]

\begin{lemma}[Positivity and monotonicity]
\label{CHlem_cross}
The renormalized Jacobian \(R\) remains strictly positive on
\(\mathcal K\). More precisely, there exist constants
\(c_R,C_R>0\) such that
\begin{equation}\label{CH_R_uniform_bounds}
    0<c_R\leq R(a,t)\leq C_R
\end{equation}
for all \((a,t)\in\mathcal K\). Consequently,
\[
    J(a,t)\geq0,
    \qquad
    M(a,t)>0,
\]
and
\begin{equation}\label{CH_J_zero_characterization}
    J(a,t)=0
    \quad\Longleftrightarrow\quad
    Q(a,t)=0.
\end{equation}
In particular, for each fixed \(t\), the map
\(a\mapsto\eta(a,t)\) is nondecreasing. If the zero set of
\(Q(\cdot,t)\) contains no nontrivial interval, then
\(\eta(\cdot,t)\) is strictly increasing on \(I_*\).
\end{lemma}

%

\begin{proof}
Equation \eqref{CH_R_Q_evolution} has the explicit representation
\begin{equation*}
    R(a,t)
    =
    R(a,0)
    \exp\left(
        \int_0^t
        \bigl(1+2\mathcal F(a,s)\bigr)
        \frac{Q(a,s)}{1+Q(a,s)^2}\,ds
    \right).
\end{equation*}
Since
\(
    R(a,0)=1+\bigl(u_0'(a)\bigr)^2>0
\)
and the coefficient in the exponential is bounded on \(\mathcal K\),
the bounds \eqref{CH_R_uniform_bounds} follow.  The formulas
\eqref{CH_J_Q_R} then imply
\[
    J=\frac{Q^2R}{1+Q^2}\geq0,
    \qquad
    M=\frac{R}{1+Q^2}>0,
\]
as well as \eqref{CH_J_zero_characterization}.

For \(a_1<a_2\),
\[
    \eta(a_2,t)-\eta(a_1,t)
    =
    \int_{a_1}^{a_2}J(a,t)\,da
    \geq0.
\]
Thus \(\eta(\cdot,t)\) is nondecreasing.  If \(Q\) does not vanish on
an entire nonempty subinterval, every interval \((a_1,a_2)\) contains
a subinterval on which \(J>0\); hence the preceding inequality is
strict.  
\end{proof}


\subsubsection{Regularity bootstrap}
\label{sec:CH-regularity-bootstrap}


By Proposition~\ref{prop:CH-local-continuation}, Lemma~\ref{CHlem_cross} and \eqref{CH_J_Q_R}, the functions $Q,R,J,N,M,\eta,U,P,\Pi$ are continuous on $\mathcal K$. Building on this, the next two lemmas establish higher regularity. 

\begin{lemma}[Regularity of the extended Lagrangian variables]
\label{lem:CH-extended-regularity}
On the compact cylinder \(\mathcal K\),
\[
    P,P_a,\mathcal F,\mathcal F_a\in C(\mathcal K),
\]
with
\begin{equation*}
    P_a(a,t)=\Pi(a,t)J(a,t)
\end{equation*}
and
\begin{equation*}
    \mathcal F_a(a,t)
    =
    2\bigl(U(a,t)+\gamma\bigr)
    \frac{Q(a,t)R(a,t)}{1+Q(a,t)^2}
    -
    \Pi(a,t)
    \frac{Q(a,t)^2R(a,t)}{1+Q(a,t)^2}.
\end{equation*}
Moreover,
\[
    Q,R,J,N,M\in C^1(\mathcal K),
    \qquad
    \eta,U\in C^1(\mathcal K),
\]
and, for each fixed \(t\),
\[
    \eta(\cdot,t),U(\cdot,t)\in C^2(I_*).
\]
In particular, for every \(a_0\in I_*\),
\begin{equation}\label{CH_R_first_order_expansion}
    R(a,t)
    =
    R(a_0,t)
    +
    O\bigl(|a-a_0|\bigr)
\end{equation}
as \(a\to a_0\), uniformly for
\(t\in[0,t_*+\varepsilon_0]\).
\end{lemma}

\begin{proof}
Differentiating \eqref{CH_extended_P} and \eqref{CH_extended_Pi} with respect to \(a\)
give
\[
    P_a
    =
    \eta_a
    \int_{\mathbb R}
    G_x\bigl(\eta(a,t)-\eta(b,t)\bigr)
    \mathcal E(b,t)\,db
    =
    J\Pi
\]
and
\[
\mathcal F_a
=
2(U+\gamma)N-\Pi J,
\]
where
\[
    \mathcal E
    =
    \bigl(U^2+2\gamma U\bigr)J+\frac12M.
\]
Thus \(P_a, \mathcal F_a\) are continuous.

Differentiating \eqref{CH_Q_evolution} with respect to \(a\) yields
\begin{equation*}
    \partial_tQ_a
    =
    -\mathcal F_aQ^2
    -
    2\mathcal FQQ_a.
\end{equation*}
Because \(Q(\cdot,0)=1/u_0'\) is \(C^1\) on \(I\),
the parameter-dependent ODE gives
\[
    Q_a\in C(\mathcal K).
\]
Equation \eqref{CH_Q_evolution} also gives \(Q_t\in C(\mathcal K)\);
hence \(Q\in C^1(\mathcal K)\).

By \eqref{CH_Q_evolution}, $R_t=BR$ where
\[
    B
    :=
    \bigl(1+2\mathcal F\bigr)
    \frac{Q}{1+Q^2}.
\]
Here,
\[
    B_a
    =
    2\mathcal F_a\frac{Q}{1+Q^2}
    +
    \bigl(1+2\mathcal F\bigr)
    \frac{1-Q^2}{(1+Q^2)^2}Q_a
\]
is continuous, and it implies that
\begin{equation*}
    \partial_tR_a
    =
    BR_a+B_aR
\end{equation*}
with \(C^1\) initial data, and consequently
\[
    R\in C^1(\mathcal K).
\]
The algebraic formula
\eqref{CH_J_Q_R} now give
\[
    J,N,M\in C^1(\mathcal K).
\]
Since \(\eta_a=J\), \(U_a=N\), \(\eta_t=U\), and \(U_t=-\Pi\), the
asserted regularity of \(\eta\) and \(U\) follows.  Finally,
\eqref{CH_R_first_order_expansion} is an immediate consequence of the
uniform boundedness of \(R_a\) on \(\mathcal K\).
\end{proof}


\begin{lemma}[Second spatial regularity bootstrap]
\label{lem:CH-second-regularity-bootstrap}
Under the assumptions of
Lemma~\ref{lem:CH-extended-regularity}, the Lagrangian force \(\Pi\) is
continuously differentiable with respect to \(a\), the forcing
\(\mathcal F\) is twice continuously differentiable with respect to
\(a\), and
\[
    \Pi_a,\mathcal F_{aa},Q_{aa}\in C(\mathcal K).
\]
\end{lemma}

\begin{proof}
We first improve the spatial regularity of the Lagrangian nonlocal
force \(\Pi\).  
Fix \(a_0\in I_*\).  Since
\[
    U(a,t)-U(a_0,t)
    =
    \int_{a_0}^{a}N(\alpha,t)\,d\alpha,
\]
and \(N_t\) is continuous, differentiation with respect to \(t\) is
justified and by $U_t=-\Pi$, 
\[
    -\Pi(a,t)+\Pi(a_0,t)
    =
    \int_{a_0}^{a}
    \left(
        \frac12M(\alpha,t)
        +
        \mathcal F(\alpha,t)J(\alpha,t)
    \right)
    \,d\alpha.
\]
Therefore,
\begin{equation}\label{CH_Pia_formula}
    \Pi_a
    =
    -\frac12M-\mathcal FJ.
\end{equation}
Since \(M\), \(\mathcal F\), and \(J\) are continuous on
\(\mathcal K\), it follows that
\[
    \Pi_a\in C(\mathcal K).
\]

Recall that
\[
    \mathcal F_a
    =
    2(U+\gamma)N-\Pi J.
\]
Differentiating with respect to \(a\), we obtain by \eqref{CH_Pia_formula} that
\begin{equation*}
	\begin{split}
		 \mathcal F_{aa}
		&=
		2N^2
		+
		2(U+\gamma)N_a
		-
		\Pi_aJ
		-
		\Pi J_a
		\\
		&=
		2N^2
		+
		2(U+\gamma)N_a
		+
		\mathcal FJ^2
		+
		\frac12MJ
		-
		\Pi J_a.
	\end{split}
\end{equation*}
Every quantity on the right-hand side is continuous on
\(\mathcal K\), and hence
\[
    \mathcal F_{aa}\in C(\mathcal K).
\]

We now differentiate \eqref{CH_Q_evolution} twice with respect to \(a\).
The first derivative satisfies
\[
    \partial_tQ_a
    =
    -\mathcal F_aQ^2
    -
    2\mathcal FQQ_a.
\]
Differentiating once more gives
\begin{equation}\label{CH_Qaa_evolution}
    \partial_tQ_{aa}
    +
    2\mathcal FQ Q_{aa}
    =
    -\mathcal F_{aa}Q^2
    -
    4\mathcal F_aQ Q_a
    -
    2\mathcal FQ_a^2.
\end{equation}
Because \(\overline{I_*}\Subset I\) and \(u_0'\) does not vanish on
\(I\), it is bounded away from zero on \(\overline{I_*}\).  Since
\(u_0\in H^5(\mathbb R)\), the initial inverse slope
\[
    Q(a,0)=\frac{1}{u_0'(a)}
\]
belongs to \(C^2(\overline{I_*})\).  Equation
\eqref{CH_Qaa_evolution} therefore yields
\[
    Q_{aa}\in C(\mathcal K).
\]
\end{proof}

\begin{corollary}[Quadratic expansion at the first breaking point]
\label{cor:CH-quadratic-Q-expansion}
Under Assumption~\ref{ass:CH-main}, one has
\begin{equation}\label{CH_Q_quadratic_expansion}
    Q(a,t_*)
    =
    -\frac14\tau''(a_*)(a-a_*)^2
    +
    o\bigl(|a-a_*|^2\bigr)
\end{equation}
as \(a\to a_*\).  Consequently,
\begin{equation}\label{CH_Qaa_tau_relation}
    Q_a(a_*,t_*)=0,
    \qquad
    Q_{aa}(a_*,t_*)
    =
    -\frac12\tau''(a_*)
    <
    0.
\end{equation}
\end{corollary}

\begin{proof}
For \(a\) near \(a_*\), the breaking-time map satisfies
\[
    Q(a,\tau(a))=0,
    \qquad
    \tau(a)
    =
    t_*
    +
    \frac12\tau''(a_*)(a-a_*)^2
    +
    o\bigl(|a-a_*|^2\bigr).
\]
Integrating \eqref{CH_Q_evolution} from \(t_*\) to \(\tau(a)\), we obtain
\[
    Q(a,t_*)
    =
    -\frac12\bigl(\tau(a)-t_*\bigr)
    +
    \int_{t_*}^{\tau(a)}
    \mathcal F(a,s)Q(a,s)^2\,ds.
\]
Since \(Q(a,\tau(a))=0\) and \(Q_t\) is bounded in a neighborhood of
\((a_*,t_*)\),
\[
    |Q(a,s)|
    \leq
    C|\tau(a)-s|
    \leq
    C|\tau(a)-t_*|
\]
for \(s\) between \(t_*\) and \(\tau(a)\).  Hence, by the boundedness of $\mathcal F$ in \eqref{CH_forcing_bound}, we obtain
\[
\left|
\int_{t_*}^{\tau(a)}
\mathcal F(a,s)Q(a,s)^2\,ds
\right|
\le C|\tau(a)-t_*|^3.
\]
Since $\tau(a)-t_*=O(|a-a_*|^2)$, it follows that
\[
\int_{t_*}^{\tau(a)}
\mathcal F(a,s)Q(a,s)^2\,ds
=o(|a-a_*|^2).
\]
Substituting the expansion of $\tau$ into the preceding identity
therefore yields \eqref{CH_Q_quadratic_expansion}.

By Lemma~\ref{lem:CH-second-regularity-bootstrap},
\(Q(\cdot,t_*)\in C^2(\overline{I_*})\).  Comparing
\eqref{CH_Q_quadratic_expansion} with the second-order Taylor formula
for \(Q(\cdot,t_*)\) at \(a_*\) yields
\eqref{CH_Qaa_tau_relation}.
\end{proof}

\begin{remark}[Strict monotonicity after the first breaking time]
\label{rmk:CH-Q-zero-set}
By the definition of the local breaking-time map,
\[
    Q(a,t)=0
    \quad\Longleftrightarrow\quad
    t=\tau(a)
\]
near \((a_*,t_*)\).
Since \(\tau''(a_*)>0\) and \(\tau\in C^2\), after shrinking \(I_*\)
if necessary, \(\tau\) is strictly convex on \(I_*\). Consequently,
for every \(t\in[t_*,t_*+\varepsilon_*)\), the zero set
\[
    \{a\in I_*:Q(a,t)=0\}
    =
    \{a\in I_*:\tau(a)=t\}
\]
contains at most two points and hence no nontrivial interval.
It follows from Lemma~\ref{CHlem_cross} that
\(\eta(\cdot,t)\) is strictly increasing on \(I_*\).
Thus, although \(J=\eta_a\) may vanish at isolated breaking labels,
no nontrivial Lagrangian interval collapses to a single Eulerian point.
\end{remark}



\subsubsection{Eulerian reconstruction}
\label{sec:CH-eulerian-reconstruction}

We now reconstruct the Eulerian velocity by means of a generalized inverse
of the Lagrangian flow map.

\begin{lemma}[Constancy on the fibers of the flow]
\label{lem:CH-constant-on-fibers}
Fix
\(
    0\leq t<t_*+\varepsilon_*.
\)
If
\(
    a_1<a_2,
    \eta(a_1,t)=\eta(a_2,t),
\)
then
\(
    U(a_1,t)=U(a_2,t).
\)
More generally, for every \(x\in\mathbb R\), \(U(\cdot,t)\) is constant
on the fiber
\[
    \eta(\cdot,t)^{-1}(\{x\})
    =
    \{a\in\mathbb R:\eta(a,t)=x\}.
\]
\end{lemma}

\begin{proof}
Since \(J=\eta_a\geq0\),
\[
    0
    =
    \eta(a_2,t)-\eta(a_1,t)
    =
    \int_{a_1}^{a_2}J(a,t)\,da.
\]
Because \(J(\cdot,t)\) is continuous and nonnegative, it follows that
\[
    J(a,t)=0
    \qquad
    \text{for all }a\in[a_1,a_2].
\]
The compatibility relation \eqref{CH_compatibility_H} then gives
\[
    N(a,t)^2=J(a,t)M(a,t)=0
\]
on \([a_1,a_2]\). Hence \(N=0\) there, and since \(U_a=N\), the function
\(U(\cdot,t)\) is constant on \([a_1,a_2]\).
\end{proof}

By Proposition~\ref{prop:CH-local-continuation}, for every fixed
\(
    0\leq t<t_*+\varepsilon_*,
\)
the map \(\eta(\cdot,t)\) is continuous and nondecreasing. Moreover,
\[
    \eta(\cdot,t)-\operatorname{Id}\in L^\infty(\mathbb R),
\]
and therefore
\[
    \eta(a,t)\longrightarrow\pm\infty
    \qquad
    \text{as }a\longrightarrow\pm\infty.
\]
Thus \(\eta(\cdot,t)\) maps \(\mathbb R\) onto \(\mathbb R\).

For \(x\in\mathbb R\), define the generalized inverse
\begin{equation}\label{CH_generalized_inverse}
    \xi(x,t)
    :=
    \sup\left\{
        a\in\mathbb R:
        \eta(a,t)\leq x
    \right\}.
\end{equation}
The set in \eqref{CH_generalized_inverse} is nonempty and bounded above, so
\(\xi(x,t)\) is well defined. By the continuity and monotonicity of
\(\eta(\cdot,t)\),
\begin{equation*}
    \eta(\xi(x,t),t)=x.
\end{equation*}

We then define the Eulerian velocity by
\begin{equation}\label{CH_eulerian_reconstruction}
    u(x,t)
    :=
    U(\xi(x,t),t).
\end{equation}
If \(a\) is any Lagrangian label satisfying
\[
    \eta(a,t)=x,
\]
then \(a\) and \(\xi(x,t)\) belong to the same fiber of
\(\eta(\cdot,t)\). Lemma~\ref{lem:CH-constant-on-fibers} therefore gives
\[
    U(a,t)=U(\xi(x,t),t)=u(x,t).
\]
Consequently,
\begin{equation*}
    U(a,t)=u(\eta(a,t),t)
\end{equation*}
for every \(a\in\mathbb R\) and
\(0\leq t<t_*+\varepsilon_*\).

For each fixed \(t\), define the nonnegative Radon measure
\[
    \mu_t
    :=
    \eta(\cdot,t)_{\#}
    \bigl(M(\cdot,t)\,da\bigr),
\]
that is, for every Borel set \(B\subset\mathbb R\),
\begin{equation*}
    \mu_t(B)
    =
    \int_{\eta(\cdot,t)^{-1}(B)}
    M(a,t)\,da.
\end{equation*}
Equivalently, for every bounded Borel function
\(\varphi:\mathbb R\to\mathbb R\),
\begin{equation*}
    \int_{\mathbb R}\varphi(x)\,d\mu_t(x)
    =
    \int_{\mathbb R}
    \varphi\bigl(\eta(a,t)\bigr)M(a,t)\,da.
\end{equation*}
Here and after, \(d\mu_t^{\mathrm{ac}}\) denotes the absolutely continuous part
of \(d\mu_t\), namely
\[
    d\mu_t^{\mathrm{ac}}
    =
    \frac{d\mu_t}{dx}\,dx. 
\]
\begin{proposition}[Eulerian reconstruction]
\label{prop:CH-eulerian-reconstruction}
For every \(T<t_*+\varepsilon_*\), the function \(u\) defined by
\eqref{CH_eulerian_reconstruction} is continuous on
\(\mathbb R\times[0,T]\). 
Then
\begin{equation}\label{CH_energy_measure_ac}
    d\mu_t^{\mathrm{ac}}
    =
    u_x(x,t)^2\,dx.
\end{equation}
The potential satisfies
\begin{equation}\label{CH_P_eulerian_measure_representation}
    P(a,t)=p(\eta(a,t),t),
\end{equation}
where
\begin{equation}\label{CH_reconstructed_p_measure}
    \begin{split}
    p(x,t)
    &:=
    \int_{\mathbb R}
    G(x-y)
    \bigl(u(y,t)^2+2\gamma u(y,t)\bigr)\,dy
    +
    \frac12
    \int_{\mathbb R}G(x-y)\,d\mu_t(y).
    \end{split}
\end{equation}
If \(\{J(\cdot,t)=0\}\) has measure zero, then
\begin{equation}\label{CH_no_concentrated_energy}
    d\mu_t=u_x(x,t)^2\,dx,
\end{equation}
and hence
\begin{equation}\label{CH_reconstructed_standard_p}
    p
    =
    (1-\partial_x^2)^{-1}
    \left(
        u^2+\frac12u_x^2+2\gamma u
    \right).
\end{equation}
Finally, the reconstructed solution agrees with the classical solution
for \(0\leq t<t_*\).
\end{proposition}

\begin{proof}
Continuity follows from the continuity of \((\eta,U)\), the boundedness
of \(\eta-\operatorname{Id}\), and the constancy of \(U\) on the fibers
of \(\eta\).

The assertions concerning \(\mu_t\) follow from the decomposition
\[
    \mu_t=\nu_t+\sigma_t,
\]
where
\[
    \nu_t
    :=
    \eta(\cdot,t)_{\#}
    \bigl(M\mathbf 1_{\{J>0\}}\,da\bigr),
    \qquad
    \sigma_t
    :=
    \eta(\cdot,t)_{\#}
    \bigl(M\mathbf 1_{\{J=0\}}\,da\bigr).
\]
On \(\{J>0\}\), we have
\[
    M(a,t)
    =
    u_x(\eta(a,t),t)^2J(a,t).
\]
Hence, by the monotone change-of-variables formula,
\[
    d\nu_t
    =
    u_x(x,t)^2\,dx.
\]

On the other hand, \(\sigma_t\) is concentrated on $\eta\bigl(\{J(\cdot,t)=0\},t\bigr)$. Since \(a\mapsto\eta(a,t)\) is absolutely continuous and
\(\eta_a=J\geq0\), the one-dimensional area formula gives
\[
    \left|
        \eta\bigl(\{J(\cdot,t)=0\},t\bigr)
    \right|
    \leq
    \int_{\{J(\cdot,t)=0\}}J(a,t)\,da
    =
    0.
\]
Thus \(\sigma_t\) is singular with respect to Lebesgue measure.
Consequently, \(\nu_t\) and \(\sigma_t\) are precisely the absolutely
continuous and singular parts of \(\mu_t\), respectively, and
\[
    d\mu_t^{\mathrm{ac}}
    =
    u_x(x,t)^2\,dx,
\]
which proves \eqref{CH_energy_measure_ac}.

If, in addition, \(\{J(\cdot,t)=0\}\) has Lebesgue measure zero, then
\[
    \sigma_t(\mathbb R)
    =
    \int_{\{J(\cdot,t)=0\}}M(a,t)\,da
    =
    0.
\]
Hence \(\sigma_t=0\), and therefore
\[
    d\mu_t=u_x(x,t)^2\,dx,
\]
which proves \eqref{CH_no_concentrated_energy}.

Finally, applying the change-of-variables formula to the definition of
\(P\) gives \eqref{CH_P_eulerian_measure_representation} and
\eqref{CH_reconstructed_p_measure}. In the absence of concentrated
energy, \eqref{CH_reconstructed_standard_p} follows immediately.

Before \(t_*\), \(\eta(\cdot,t)\) is a diffeomorphism and the generalized
inverse is the ordinary inverse. Therefore, the reconstruction agrees
with the classical solution.
\end{proof} 

\begin{remark}
\label{rmk:CH-reconstruction-before-geometry}
Under Assumption~\ref{ass:CH-main}, the generalized-inverse
construction does not require 
$\{J(\cdot,t)=0\}$
to be known a priori to have measure zero. The cusp-pair analysis below
shows that the relevant zeros are isolated, and hence that the
generalized inverse coincides locally with the ordinary inverse.
\end{remark}

\subsection{Formation of a local anticuspon--cuspon pair}
\label{sec_CHcusp}

\begin{lemma}[Splitting of the first breaking label]
	\label{CHlem:cusp-pair-no-theta}
	Let \(Q\) be the extended Lagrangian variable determined by
	\eqref{CH_Q_evolution}, and let \((a_*,t_*)\) denote the first breaking
	label and time. Assume that the first breaking is nondegenerate in the
	sense of Assumption~\ref{ass:CH-main}.
	After decreasing \(\varepsilon_*>0\), if necessary, there exist \(\delta>0\) and two functions
	\[
	a_\pm(t)
	\in
	C^1\bigl((t_*,t_*+\varepsilon_*)\bigr)
	\]
	such that
	\[
	a_-(t)<a_*<a_+(t),
	\qquad
	Q(a_\pm(t),t)=0,
	\qquad
	t\in(t_*,t_*+\varepsilon_*).
	\]
	Moreover,
	\[
	\lim_{t\searrow t_*}a_\pm(t)=a_*,
	\]
	and the two zeros are simple, with
	\[
	Q_a(a_-(t),t)>0,
	\qquad
	Q_a(a_+(t),t)<0.
	\]

	For every \(t\in(t_*,t_*+\varepsilon_*)\), one has
	\[
	Q(a,t)>0
	\qquad
	\text{for }
	a\in(a_-(t),a_+(t)),
	\]
	whereas
	\[
	Q(a,t)<0\qquad
	\text{for }a
	\in
	\bigl(a_-(t)-\delta,a_-(t)\bigr)
	\cup
	\bigl(a_+(t),a_+(t)+\delta\bigr).
	\]

	Finally, as \(t\searrow t_*\), the two breaking labels separate according
	to
	\begin{equation}\label{apm_ext}
		a_\pm(t)
		=
		a_*
		\pm
		\sqrt{
			\frac{t-t_*}
			{-Q_{aa}(a_*,t_*)}
		}
		+
		o\bigl((t-t_*)^{1/2}\bigr).
	\end{equation}
\end{lemma}

\begin{proof}
	By Assumption~\ref{ass:CH-main}, there exists a \(C^2\)
	breaking-time map \(\tau\), defined in a neighborhood of \(a_*\), such
	that
	\[
	Q(a,\tau(a))=0,
	\qquad
	\tau(a_*)=t_*,
	\qquad
	\tau'(a_*)=0,
	\qquad
	\tau''(a_*)>0.
	\]
	Corollary~\ref{cor:CH-quadratic-Q-expansion} gives
	\begin{equation*}
		Q(a_*,t_*)=0,
		\qquad
		Q_a(a_*,t_*)=0,
		\qquad
		Q_{aa}(a_*,t_*)
		=
		-\frac12\tau''(a_*)
		<
		0.
	\end{equation*}
	In particular, \(Q(\cdot,t_*)\) has a nondegenerate local maximum at
	\(a_*\).  Moreover, equation~\eqref{CH_Q_evolution} gives
	\begin{equation*}
		Q_t(a_*,t_*)
		=
		\frac12.
	\end{equation*}

	We next construct the two post-breaking zeros of $Q$. Since
	\(\tau''(a_*)>0\), after shrinking the neighborhood of \(a_*\), if
	necessary, there exists \(\delta>0\) such that
	\begin{subequations}\label{tau_sign}
		\begin{align}
			&\tau'(a)<0
			\quad\text{for }
			a\in(a_*-\delta,a_*),
			\\
			&\tau'(a)>0
			\quad\text{for }
			a\in(a_*,a_*+\delta).
		\end{align}
	\end{subequations}
	Hence \(\tau\) is strictly decreasing to the left of \(a_*\) and
	strictly increasing to the right of \(a_*\). Since
	\(\tau(a_*)=t_*\), after decreasing \(\varepsilon_*\), if necessary, for each
	\(t\in(t_*,t_*+\varepsilon_*)\), the equation
	\(
	\tau(a)=t
	\) 
	has exactly two solutions in
	\((a_*-\delta,a_*+\delta )\), denoted by
	\(
	a_-(t)<a_*<a_+(t).
	\) 
	Because the local zero set of \(Q\) is the graph \(t=\tau(a)\), these
	solutions satisfy
	\(
	Q(a_\pm(t),t)=0.
	\) 
	Moreover,
	\[
	\lim_{t\searrow t_*}a_\pm(t)=a_*.
	\]
	Since
	\(
	\tau'(a_\pm(t))\neq0,
	\)
	the inverse function theorem applied separately to the two monotone
	branches of \(\tau\) shows that
	\[
	a_\pm
	\in
	C^1\bigl((t_*,t_*+\varepsilon_*)\bigr).
	\]

	We now verify that the two zeros are simple.  Differentiating
	\(
	    Q(a,\tau(a))=0
	\) 
	once with respect to \(a\) gives
	\[
	    Q_a(a,\tau(a))
	    +
	    Q_t(a,\tau(a))\tau'(a)
	    =
	    0.
	\]
	At any point of the breaking-time curve,
	equation~\eqref{CH_Q_evolution} gives
	\[
	Q_t(a,\tau(a))
	=
	\frac12,
	\]
	because \(Q(a,\tau(a))=0\). Therefore,
	\begin{equation}\label{CH_Qa_tau}
		Q_a(a,\tau(a))
		=
		-\frac12\tau'(a).
	\end{equation}
	By \eqref{tau_sign}, we have
	\[
	\tau'(a_-(t))<0,
	\qquad
	\tau'(a_+(t))>0.
	\]
	It follows from \eqref{CH_Qa_tau} that
	\[
	Q_a(a_-(t),t)>0,
	\qquad
	Q_a(a_+(t),t)<0.
	\]
	In particular, both post-breaking zeros are simple.

	We next determine the sign of \(Q\). By continuity, after decreasing
	\(\delta \) and \(\varepsilon_*\), if necessary, equation
	\eqref{CH_Q_evolution} implies
	\(
	Q_t(a,s)>0
	\)
	in a sufficiently small neighborhood of the local zero set
	\(s=\tau(a)\). Hence
	\[
	Q(a,t)
	=
	\int_{\tau(a)}^t Q_t(a,s)\,ds
	\]
	has the same sign as \(t-\tau(a)\).
	For fixed \(t\in(t_*,t_*+\varepsilon_*)\), the monotonicity of \(\tau\)
	gives
	\[
	\tau(a)<t
	\quad\text{for }
	a\in(a_-(t),a_+(t)),
	\]
	and therefore
	\[
	Q(a,t)>0
	\quad\text{for }
	a\in(a_-(t),a_+(t)).
	\]
	Similarly,
	\[
	Q(a,t)<0\quad\text{for }a\in
	\bigl(a_-(t)-\delta,a_-(t)\bigr)
	\cup
	\bigl(a_+(t),a_+(t)+\delta\bigr).
	\]

	Finally, Taylor expansion of the breaking-time map at \(a_*\) gives
	\[
	\tau(a)
	=
	t_*
	+
	\frac12\tau''(a_*)(a-a_*)^2
	+
	o\bigl(|a-a_*|^2\bigr).
	\]
	Substituting \(a=a_\pm(t)\) and using
	\(\tau(a_\pm(t))=t\), we obtain
	\[
	t-t_*
	=
	\frac12\tau''(a_*)
	\bigl(a_\pm(t)-a_*\bigr)^2
	+
	o\bigl(|a_\pm(t)-a_*|^2\bigr).
	\]
	Since
	\( \textstyle 
	Q_{aa}(a_*,t_*)
	=
	-\frac12\tau''(a_*),
	\) 
	it follows that
	\[
	\bigl(a_\pm(t)-a_*\bigr)^2
	=
	\frac{t-t_*}{-Q_{aa}(a_*,t_*)}
	+
	o(t-t_*).
	\]
	The inequalities
	\(
	a_-(t)<a_*<a_+(t)
	\)
	then determine the two signs and yield
	\[
	a_\pm(t)
	=
	a_*
	\pm
	\sqrt{
		\frac{t-t_*}
		{-Q_{aa}(a_*,t_*)}
	}
	+
	o\bigl((t-t_*)^{1/2}\bigr),
	\]
	which is \eqref{apm_ext}.
	This completes the proof. 
\end{proof}
Set 
\begin{equation}\label{x-pm}
    x_\pm(t):=\eta(a_\pm(t),t).
\end{equation} 
\begin{corollary}[Local Eulerian anticuspon--cuspon geometry]
\label{cor:CH-local-cusp-geometry}
After decreasing \(\varepsilon_*>0\), if necessary, let
\(
    I_t
    :=
    \bigl(a_-(t)-\delta,\ a_+(t)+\delta\bigr).
\)
For every \(t\in(t_*,t_*+\varepsilon_*)\), the map
\(
    \eta(\cdot,t):I_t\longrightarrow\eta(I_t,t)
\)
is strictly increasing.  Hence the generalized inverse in
\eqref{CH_generalized_inverse} coincides on \(\eta(I_t,t)\) with the
ordinary monotone inverse.

Then \(x_-(t)<x_+(t)\), and the local singular set of the reconstructed
Eulerian profile  is precisely
\(
    \{x_-(t),x_+(t)\}.
\)
The one-sided slopes satisfy
\[
    \lim_{x\to x_-(t)^-}u_x(x,t)=-\infty,
    \qquad
    \lim_{x\to x_-(t)^+}u_x(x,t)=+\infty,
\]
and
\[
    \lim_{x\to x_+(t)^-}u_x(x,t)=+\infty,
    \qquad
    \lim_{x\to x_+(t)^+}u_x(x,t)=-\infty.
\]
Moreover, \(x_-(t)\) is a cusped local minimum and \(x_+(t)\) is a
cusped local maximum.  Thus the left singularity is an anticuspon and
the right singularity is a cuspon.
\end{corollary}

\begin{proof}
By Lemma~\ref{CHlem_cross},
\(\textstyle 
    J
    =
    \frac{Q^2R}{1+Q^2}
    \geq0,
\) 
and \(J=0\) precisely where \(Q=0\).  Lemma
\ref{CHlem:cusp-pair-no-theta} shows that the only zeros of \(Q\) in
\(I_t\) are the two isolated labels \(a_-(t)\) and \(a_+(t)\).
Therefore, for \(a_1<a_2\) in \(I_t\),
\[
    \eta(a_2,t)-\eta(a_1,t)
    =
    \int_{a_1}^{a_2}J(a,t)\,da
    >
    0.
\]
Thus \(\eta(\cdot,t)\) is strictly increasing on \(I_t\).

On the regular set, Proposition~\ref{prop:CH-eulerian-reconstruction}
gives
\( \textstyle 
    u_x(\eta(a,t),t)=\frac1{Q(a,t)}.
\) 
The sign structure in Lemma~\ref{CHlem:cusp-pair-no-theta} and the
simplicity of the two zeros yield the four one-sided limits above.
Away from the two zeros, \(Q\neq0\), and the Eulerian profile is
classical; hence no other local singularities occur.

Finally,
\( \textstyle 
    U_a
    =
    \frac{QR}{1+Q^2}
\) 
has the same sign as \(Q\).  Thus \(U(\cdot,t)\) decreases immediately
to the left of \(a_-(t)\), increases on
\((a_-(t),a_+(t))\), and decreases immediately to the right of
\(a_+(t)\).  Since \(\eta(\cdot,t)\) is increasing, the same
monotonicity pattern holds for the Eulerian profile.
\end{proof}

\begin{remark}[Local monotonicity reversal]
\label{rmk:CH-local-monotonicity-reversal}
For \(t<t_*\) sufficiently close to \(t_*\), the local normal form gives
\( \textstyle 
    Q(a,t)<0
\) 
near \(a_*\), and hence \(U_a<0\).  The profile is therefore locally
decreasing near the central characteristic before breaking.  For
\(t>t_*\), Lemma~\ref{CHlem:cusp-pair-no-theta} implies that
\[
    Q(a,t)>0
    \qquad
    \text{for }a\in(a_-(t),a_+(t)),
\]
so that the portion of the profile between the newly formed singular
points is strictly increasing.  Thus the anticuspon--cuspon
bifurcation is accompanied by an immediate local reversal of
monotonicity.
\end{remark}


\subsection{Sharp H\"older regularity and outward motion of the
anticuspon--cuspon pair}\label{sec_CH_reg}
We now establish the sharp H\"older regularity of the anticuspon-cuspon pair. 
Let \( x_* := \eta(a_*, t_*)$ and $x_\pm(t)$ be as defined in \eqref{x-pm}. 
\begin{lemma}[Sharp pointwise H\"older exponents at the cusp points]
\label{CHreg_lem}
Under the assumptions of
Lemma~\ref{CHlem:cusp-pair-no-theta}, there exists a constant
\(C_*\in(0,\infty)\) such that
\[
    \lim_{a\to a_*}
    \frac{
        \bigl|
            u(\eta(a,t_*),t_*)
            -
            u(x_*,t_*)
        \bigr|
    }{
        \bigl|
            \eta(a,t_*)-x_*
        \bigr|^{3/5}
    }
    =
    C_*.
\]
Moreover, for every fixed \(t\in(t_*,t_*+\varepsilon_*)\) and each
\(\sigma\in\{-,+\}\), there exists
\(C_\sigma(t)\in(0,\infty)\) such that
\[
    \lim_{a\to a_\sigma(t)}
    \frac{
        \bigl|
            u(\eta(a,t),t)
            -
            u(x_\sigma(t),t)
        \bigr|
    }{
        \bigl|
            \eta(a,t)-x_\sigma(t)
        \bigr|^{2/3}
    }
    =
    C_\sigma(t).
\]
Consequently, the first breaking profile has exact pointwise H\"older
exponent \(3/5\) at \(x_*\), while each post-breaking cusp has exact
pointwise H\"older exponent \(2/3\).
\end{lemma}

\begin{proof}
Let \((a_0,t_0)\) be
either
\(
    (a_0,t_0)=(a_*,t_*),
\)
in which case
\[
    Q(a_*,t_*)=Q_a(a_*,t_*)=0,
    \qquad
    Q_{aa}(a_*,t_*)\neq0,
\]
or
\(
    (a_0,t_0)=(a_\sigma(t),t),
    \sigma\in\{-,+\},
\)
in which case
\[
    Q(a_\sigma(t),t)=0,
    \qquad
    Q_a(a_\sigma(t),t)\neq0.
\]
Accordingly, as \(h\to0\),
\begin{equation}\label{mq}
	  Q(a_0+h,t_0)
	=
	qh^m+o\bigl(|h|^m\bigr),
	\qquad
	q\neq0,
\end{equation}
where
\[
    (m,q)
    =
    \begin{cases}
        \displaystyle
        \left(
            2,\frac12Q_{aa}(a_*,t_*)
        \right),
        & (a_0,t_0)=(a_*,t_*),\\[1.2ex]
        \displaystyle
        \left(
            1,Q_a(a_\sigma(t),t)
        \right),
        & (a_0,t_0)=(a_\sigma(t),t).
    \end{cases}
\]

By
Lemma~\ref{lem:CH-extended-regularity} and $R>0$ by Lemma~\ref{CHlem_cross},
\[
    R(a_0+h,t_0)
    =
    r+O(|h|),
    \qquad
    r:=R(a_0,t_0)>0.
\]
Since for \(|Q|\ll1\),
\[
    \frac{Q}{1+Q^2}
    =
    Q+O(|Q|^3),
    \qquad
    \frac{Q^2}{1+Q^2}
    =
    Q^2+O(Q^4),
\]
we obtain by \eqref{CH_J_Q_R} that
\[
    U_a(a_0+h,t_0)
    =
    rqh^m+o\bigl(|h|^m\bigr)
\]
and
\[
    \eta_a(a_0+h,t_0)
    =
    rq^2h^{2m}+o\bigl(|h|^{2m}\bigr).
\]
Integrating from \(0\) to \(h\), we find
\begin{equation*}
    U(a_0+h,t_0)-U(a_0,t_0)
    =
    \frac{rq}{m+1}h^{m+1}
    +
    o\bigl(|h|^{m+1}\bigr)
\end{equation*}
and
\begin{equation*}
    \eta(a_0+h,t_0)-\eta(a_0,t_0)
    =
    \frac{rq^2}{2m+1}h^{2m+1}
    +
    o\bigl(|h|^{2m+1}\bigr).
\end{equation*}
It follows that
\begin{equation}\label{CH_holder_local_limit}
    \lim_{h\to0}
    \frac{
        \bigl|
            U(a_0+h,t_0)-U(a_0,t_0)
        \bigr|
    }{
        \bigl|
            \eta(a_0+h,t_0)-\eta(a_0,t_0)
        \bigr|^{(m+1)/(2m+1)}
    }
    =
    \frac{
        r|q|/(m+1)
    }{
        \bigl(rq^2/(2m+1)\bigr)^{(m+1)/(2m+1)}
    }
    \in(0,\infty).
\end{equation}
Combining \eqref{mq} and \eqref{CH_holder_local_limit} completes the proof.
%
\end{proof}


The sign information for \(Q\) obtained in the construction of the
anticuspon--cuspon pair also determines the direction of motion of the
two singular branches.  In particular, we obtain the following
geometric consequence.

\begin{corollary}[Outward motion of the anticuspon--cuspon pair]
\label{cor:CH-outward-motion}
Let
\[
    x_\pm(t)=\eta(a_\pm(t),t),
    \qquad
    x_c(t)=\eta(a_*,t).
\]
Then, for \(t\in(t_*,t_*+\varepsilon_*)\),
\[
    \frac{d}{dt}\bigl(x_+(t)-x_c(t)\bigr)>0,
    \qquad
    \frac{d}{dt}\bigl(x_c(t)-x_-(t)\bigr)>0.
\]
Hence the anticuspon--cuspon pair separates strictly outward from the
central characteristic immediately after breaking.
\end{corollary}

\begin{proof}
Since \(J(a_\pm(t),t)=0\) and \(\eta_t=U\),
\[
    x_\pm'(t)=U(a_\pm(t),t),
    \qquad
    x_c'(t)=U(a_*,t).
\]
Using \eqref{CH_J_Q_R},
\[
\begin{aligned}
x_+'(t)-x_c'(t)
&=
\int_{a_*}^{a_+(t)}
\frac{QR}{1+Q^2}\,da>0,\\
x_c'(t)-x_-'(t)
&=
\int_{a_-(t)}^{a_*}
\frac{QR}{1+Q^2}\,da>0,
\end{aligned}
\]
because \(Q>0\) on \((a_-(t),a_+(t))\) and \(R>0\).
\end{proof}

\subsection{Identification with the conservative solution}
\label{sec_u_unique}

We now show that the Eulerian solution reconstructed from the extended
Lagrangian variables belongs to the conservative solution class of
\cite{BCZ}. Throughout this subsection, let
\(
    0<T<t_*+\varepsilon_*.
\) 
To formulate the conservation of the derivative energy across the cusp
curves, set
\(
    w:=u_x^2.
\)
For smooth solutions, differentiating the equation with respect to \(x\)
gives the balance law
\begin{equation}\label{balance}
    w_t+(uw)_x
    =
    2\bigl(u^2+2\gamma u-p\bigr)u_x.
\end{equation}

\begin{definition}[Conservative solution]
\label{CH_cons_def}
A solution \(u\) of
\eqref{CH_init}--\eqref{CH_nonlocal} on \([0,T]\) is called conservative
if \(w=u_x^2\) satisfies \eqref{balance} in the sense of distributions.
Equivalently, for every
\(
    \psi\in C_c^1(\mathbb{R}\times [0,T)),
\) 
one has
\[
    \int_0^T\int_{\mathbb R}
    \left[
        u_x^2\psi_t
        +
        uu_x^2\psi_x
        +
        2\bigl(u^2+2\gamma u-p\bigr)u_x\psi
    \right]
    \,dx\,dt
    +
    \int_{\mathbb R}
        u_{0,x}^2(x)\psi(0,x)\,dx
    =
    0.
\]
\end{definition}

We shall also use the solution class introduced in \cite{BCZ}.

\begin{definition}
\label{CH_sol_def}
A function \(u=u(x,t)\) is a solution of
\eqref{CH_init}--\eqref{CH_nonlocal} on \([0,T]\) if the following
properties hold:
\begin{enumerate}
    \item
    \(u\) is H\"older continuous on \(\mathbb R\times [0,T]\);

    \item
    \(u(\cdot,t)\in H^1(\mathbb R)\) for every \(t\in[0,T]\);

    \item
    the map
    \(
        t\longmapsto u(\cdot,t)
    \) 
    is Lipschitz continuous from \([0,T]\) into \(L^2(\mathbb R)\);

    \item
    \(u(\cdot,0)=u_0\), and
\(
    \textstyle 
        \frac{d}{dt}u
        =
        -uu_x-p_x
        \)
    in \(L^2(\mathbb R)\) for almost every \(t\in[0,T]\).
\end{enumerate}
\end{definition}

\begin{proposition}[\cite{BCZ}]
\label{CH_unique}
Let \(\gamma=0\). For every \(u_0\in H^1(\mathbb R)\), the Cauchy problem
\eqref{CH_init}--\eqref{CH_nonlocal} has a unique solution satisfying
Definitions~\ref{CH_cons_def} and \ref{CH_sol_def}.
\end{proposition}

By Assumption~\ref{ass:CH-main}, the first breaking label is unique. After decreasing \(\varepsilon_*\), if necessary, Corollary~\ref{cor:CH-local-cusp-geometry} then implies that
%
\(
    \{a\in\mathbb R:J(a,t)=0\}
\) 
has Lebesgue measure zero for every
\(t\in[0,t_*+\varepsilon_*)\).  Indeed, the flow remains classical away
from the first-breaking neighborhood, while inside that neighborhood
the zero set consists of \(a_*\) at \(t=t_*\) and of
\(a_-(t),a_+(t)\) for \(t>t_*\).  Hence
\[
    d\mu_t=u_x(x,t)^2\,dx
\]
by Proposition~\ref{prop:CH-eulerian-reconstruction}.

We first verify 
the derivative-energy balance
law, and conservation of the full \(H^1\)-energy.

\begin{lemma}[Eulerian equation and conservativity]
\label{lem_cons_sol}
Let \(u\) be the Eulerian solution reconstructed from the extended
Lagrangian variables by
\eqref{CH_eulerian_reconstruction}. Then 
\(u\) is conservative in the sense of
Definition~\ref{CH_cons_def}. In addition, it holds that 
\begin{equation}\label{CH_reconstructed_H1_conservation}
    \|u(\cdot,t)\|_{H^1(\mathbb R)}
    =
    \|u_0\|_{H^1(\mathbb R)},
    \qquad
    t\in[0,T].
\end{equation}
\end{lemma}

\begin{proof}

We first prove the balance law. Let
\(
    \psi\in C_c^1(\mathbb R\times [0,T))
\) 
and let \(\eta\) be a characteristic curve associated with \(u\). 
By \eqref{CH_derived_M}, we have 
\begin{equation*}
\begin{split}
  &  \frac{d}{dt}
    \left[
        M(a,t)\psi(\eta(a,t),t)
    \right]
\\  
& \qquad =
    M(a,t)
    \left[
        \psi_t(\eta(a,t),t)
        +
        U(a,t)\psi_x(\eta(a,t),t)
    \right]
    +
    2\mathcal F(a,t)N(a,t)
        \psi(\eta(a,t),t).
\end{split}
\end{equation*} 
Integrating in \(a\) and \(t\), and using
\(
    M(a,0)=u_{0,x}^2(a), 
    \eta(a,0)=a,
\)
we obtain
\[
    \begin{aligned}
    \int_0^T\int_{\mathbb R}
    \left[
        M(a,t)
        \frac{d}{dt}\psi(\eta(a,t),t)
        +
        2\mathcal F(a,t)N(a,t)
        \psi(\eta(a,t),t)
    \right]
    \,da\,dt
    +
    \int_{\mathbb R}
        M(a,0)\psi(a,0)\,da = 0.
    \end{aligned}
\]
Rewriting the preceding identity in Eulerian variables via the monotone change of variables $x=\eta(a,t)$, we have
\[
    \int_0^T\int_{\mathbb R}
    \left[
        u_x^2\psi_t
        +
        uu_x^2\psi_x
        +
        2\bigl(u^2+2\gamma u-p\bigr)u_x\psi
    \right]
    \,dx\,dt
    +
    \int_{\mathbb R}
        u_{0,x}^2(x)\psi(x,0)\,dx=0.
\]
Thus \(u\) is conservative in the sense of
Definition~\ref{CH_cons_def}.

Next, we prove conservation of the full \(H^1\)-energy directly in
Lagrangian variables. Since
\(
    P_a
    =
    \Pi J,
    U_a=N,
\)
we compute by \eqref{CH_regular_U} and \eqref{CH_derived_J} that 
\[
    \begin{aligned}
    \partial_t\bigl(U^2J+M\bigr)
    &=
    2UU_tJ+U^2J_t+M_t
    \\
    &=
    -2U\Pi J
    +
    U^2N
    +
    2\bigl(U^2+2\gamma U-P\bigr)N
    \\
    &=
    \partial_a
    \left(
        -2UP
        +
        U^3
        +
        2\gamma U^2
    \right).
    \end{aligned}
\]
Hence, it holds that
\[
    \frac{d}{dt}
    \int_{\mathbb R}
    \bigl(U(a,t)^2J(a,t)+M(a,t)\bigr)\,da
    =
    0.
\]
Using the monotone change of variables and the absence of a
singular part in \(\mu_t\), we have
\[
    \int_{\mathbb R}U(a,t)^2J(a,t)\,da
    =
    \int_{\mathbb R}u(x,t)^2\,dx
\]
and
\[
    \int_{\mathbb R}M(a,t)\,da
    =
    \int_{\mathbb R}u_x(x,t)^2\,dx.
\]
Consequently,
\[
    \|u(\cdot,t)\|_{H^1(\mathbb R)}^2
    =
    \|u_0\|_{H^1(\mathbb R)}^2,
\]
which proves \eqref{CH_reconstructed_H1_conservation}.
\end{proof}


We next verify the additional regularity required in
Definition~\ref{CH_sol_def}.

\begin{lemma}[Time regularity]
\label{lem_lip_t}
Let \(u\) be the reconstructed Eulerian solution. Then
\(
    u
    \in
    W^{1,\infty}(0,T;L^2(\mathbb R)),
\) 
and
\(
    u_t
    =
    -uu_x-p_x
\)
in \(L^2(\mathbb R)\) for almost every \(t\in[0,T]\). In particular, the
map
\(
    t\longmapsto u(\cdot,t)
\)
is Lipschitz continuous from \([0,T]\) into \(L^2(\mathbb R)\).
Moreover, \(u\) is H\"older continuous on
\([0,T]\times\mathbb R\). Consequently, \(u\) satisfies
Definition~\ref{CH_sol_def}.
\end{lemma}

\begin{proof}
By Lemma~\ref{lem_cons_sol},
\(\textstyle 
    \sup_{t\in[0,T]}
    \|u(\cdot,t)\|_{H^1(\mathbb R)}
    =
    \|u_0\|_{H^1(\mathbb R)}.
\)
The one-dimensional Sobolev inequality gives
\(
    \|u(\cdot,t)\|_{L^\infty}
    \leq
    C\|u_0\|_{H^1}.
\)
Therefore,
\[
    \|uu_x\|_{L^2}
    \leq
    \|u\|_{L^\infty}\|u_x\|_{L^2}
    \leq
    C\|u_0\|_{H^1}^2.
\]
Moreover, by Lemma~\ref{CH_lem_up},
\[
    \|p_x\|_{L^2}
    \leq
    C
    \left(
        \|u\|_{H^1}^2
        +
        |\gamma|\|u\|_{H^1}
    \right).
\]
Thus
\[
    f:=-uu_x-p_x
    \in
    L^\infty(0,T;L^2(\mathbb R)).
\]

By 
   Proposition~\ref{prop:CH-eulerian-reconstruction} under Assumption~\ref{ass:CH-main}, 
one has
\(
    \partial_tu=f
\) 
in the sense of distributions. Hence
\(
    u\in W^{1,\infty}(0,T;L^2(\mathbb R)),
\)
and
\(
    u_t=-uu_x-p_x
\)
in \(L^2(\mathbb R)\) for almost every \(t\in[0,T]\). In particular,
for \(0\leq s<t\leq T\),
\[
    \|u(\cdot,t)-u(\cdot,s)\|_{L^2}
    \leq
    \int_s^t
        \|f(\cdot,\tau)\|_{L^2}
    \,d\tau
    \leq
    C_{\gamma,u_0}|t-s|.
\]

It remains to verify the H\"older continuity required in
Definition~\ref{CH_sol_def}. For every fixed \(t\),
\[
    |u(x,t)-u(y,t)|
    \leq
    \|u_x(\cdot,t)\|_{L^2}
    |x-y|^{1/2}
    \leq
    C_{u_0}|x-y|^{1/2}.
\]
On the other hand, 
the one-dimensional Gagliardo--Nirenberg inequality gives
\[
    \begin{aligned}
    \|u(\cdot,t)-u(\cdot,s)\|_{L^\infty}
    &\leq
    C
    \|u(\cdot,t)-u(\cdot,s)\|_{L^2}^{1/2}
    \|\partial_xu(\cdot,t)-\partial_xu(\cdot,s)\|_{L^2}^{1/2}
    \\
    &\leq
    C_{\gamma,u_0}|t-s|^{1/2},
    \end{aligned}
\]
where we have used \(
    \|\partial_xu(\cdot,t)-\partial_xu(\cdot,s)\|_{L^2}
    \leq
    2\|u_0\|_{H^1}.
\)
Consequently,
\[
    |u(x,t)-u(y,s)|
    \leq
    C_{\gamma,u_0}
    \left(
        |x-y|^{1/2}
        +
        |t-s|^{1/2}
    \right).
\]
Thus \(u\) is H\"older continuous on
\([0,T]\times\mathbb R\). The initial condition follows directly from
the Lagrangian construction. Hence all the conditions in
Definition~\ref{CH_sol_def} are satisfied.
\end{proof}

Combining Lemmas~\ref{lem_cons_sol} and \ref{lem_lip_t}, we conclude that
the reconstructed Eulerian solution satisfies
Definitions~\ref{CH_cons_def} and \ref{CH_sol_def} on every interval
\([0,T]\) with
\(
    T<t_*+\varepsilon_*.
\)
In particular, when \(\gamma=0\), Proposition~\ref{CH_unique} implies that
the reconstructed solution coincides with the unique conservative
solution in the solution class of \cite{BCZ} throughout the entire
interval on which the reconstruction is defined.

\subsection{Completion of the proof of Theorem~\ref{Thm_CH}}
\label{sec:CH-completion}

Proposition~\ref{prop:CH-local-continuation} and
Proposition~\ref{prop:CH-eulerian-reconstruction} provide the
post-breaking continuation and the Eulerian reconstruction.
Lemma~\ref{CHlem:cusp-pair-no-theta} and
Corollary~\ref{cor:CH-local-cusp-geometry} identify the two
post-breaking singular branches and their anticuspon--cuspon
orientations.  Lemma~\ref{CHreg_lem} gives the sharp exponents
\(3/5\) at the first breaking point and \(2/3\) along the two
post-breaking branches, while Corollary~\ref{cor:CH-outward-motion} proves their outward
motion.  Finally, Lemmas~\ref{lem_cons_sol} and \ref{lem_lip_t} identify
the reconstructed profile with a conservative weak solution and, when
\(\gamma=0\), with the unique conservative solution in the class of
\cite{BCZ}.  This completes the proof of
Theorem~\ref{Thm_CH}.

\section{Generalized Hunter-Saxton Equation}\label{sec2}
This section is devoted to the proof of Theorem~\ref{mainthm1}.  
We consider the generalized Hunter--Saxton equation on the real line \(\mathbb R\):
\begin{equation}\label{genHS}
    v_{tx}+v v_{xx}
    =
    -\frac{1}{\lambda}v_x^2, \qquad \lambda \in 2\mathbb{N}. 
\end{equation}
For \(\lambda=2\), \eqref{genHS} reduces to the standard
Hunter--Saxton equation in its \(x\)-differentiated form. 
Throughout this section, we assume Assumption~\ref{ass:HS-main}.

\subsection{Post-breaking Lagrangian continuation}
\label{sec_HS_ext}

We begin with the characteristic representation of the classical solution.
Let \(\zeta(a,t)\) denote the characteristic map, defined as long as the
classical solution exists, by
\begin{equation}\label{eta_eq}
    \partial_t\zeta(a,t)
    =
    v(\zeta(a,t),t),
    \qquad
    \zeta(a,0)=a.
\end{equation}
Along each characteristic, as long as \(v_x(\zeta(a,t),t)\) remains finite,
equation~\eqref{genHS} yields
\begin{equation*}
    \frac{d}{dt}v_x(\zeta(a,t),t)
    =
    -\frac{1}{\lambda}
    v_x(\zeta(a,t),t)^2.
\end{equation*}
Integrating this Riccati equation gives
\begin{equation}\label{ux_ricc}
    v_x(\zeta(a,t),t)
    =
    \frac{v_0'(a)}
    {1+\frac{t}{\lambda}v_0'(a)}.
\end{equation}

Recall that
\begin{equation*}
    D_{\mathrm{HS}}(a,t)
    =
    1+\frac{t}{\lambda}v_0'(a).
\end{equation*}
Differentiating \eqref{eta_eq} with respect to \(a\) and using
\eqref{ux_ricc}, we obtain
\begin{equation}\label{HS_eta_a}
    \begin{aligned}
    \zeta_a(a,t)
    &=
    \exp\left(
        \int_0^t
        v_x(\zeta(a,s),s)\,ds
    \right)
   =
    \left(
        1+\frac{t}{\lambda}v_0'(a)
    \right)^\lambda
    =
    D_{\mathrm{HS}}(a,t)^\lambda.
    \end{aligned}
\end{equation}

Since \(v_0'\) attains its global minimum at \(a=a_0\),
\(D_{\mathrm{HS}}\) first vanishes at \((a_0,t_*)\), where
\begin{equation}\label{HS_tstar}
    t_*
    =
    -\frac{\lambda}{v_0'(a_0)}.
\end{equation}
In view of \eqref{ux_ricc}, \(t_*\) is precisely the first singular
time.

The explicit formulas \eqref{ux_ricc} and \eqref{HS_eta_a} identify the
Lagrangian locations at which the classical description becomes singular.
Accordingly, we define the Lagrangian singular set by
\begin{equation}\label{lagrangian_singular_set}
    \mathcal S_L
    :=
    \left\{
        (a,t)\in\mathbb R\times[0,\infty)
        :
        D_{\mathrm{HS}}(a,t)=0
    \right\}.
\end{equation}
A precise description of \(\mathcal S_L\) for times immediately after
\(t_*\) will be given in Lemma~\ref{lem:singular_set} below. 
Once the characteristic map has been continued beyond the singular time,
we denote the corresponding Eulerian singular set by
\begin{equation*}
    \mathcal S_E
    :=
    \left\{
        \bigl(\zeta(a,t),t\bigr)
        :
        (a,t)\in\mathcal S_L
    \right\}.
\end{equation*}

The equation \(D_{\mathrm{HS}}(a,t)=0\) is equivalent, whenever
\(v_0'(a)<0\), to
\[
    t=t_s(a),
    \qquad
    t_s(a):=-\frac{\lambda}{v_0'(a)}.
\]
Thus, wherever \(t_s\) is well defined and smooth, the Lagrangian
singular set is locally represented as the graph \(t=t_s(a)\).
Its Eulerian image will be considered after the Lagrangian map has
been continued beyond \(t_*\).

We now extend the Lagrangian variables beyond the first singular time.
Since \(\lambda\in2\mathbb N\), the identity
\eqref{HS_eta_a} suggests defining the extended Lagrangian map by
\begin{equation*}
    \zeta(a,t)
    :=
    \zeta(a_0,t)
    +
    \int_{a_0}^a
    D_{\mathrm{HS}}(y,t)^\lambda\,dy,
    \qquad t\geq0.
\end{equation*}
Here \(\zeta(a_0,t)\) is fixed by the corresponding Lagrangian
trajectory at the reference label \(a_0\).
On the classical lifespan, \eqref{eta_def} agrees with the original
characteristic map, while for \(t\geq t_*\) it provides its natural
continuation through the zero set of \(D_{\mathrm{HS}}\).

Similarly, on the classical lifespan, let
\[
    V_{\mathrm{HS}}(a,t)
    :=
    v(\zeta(a,t),t).
\]
Combining \eqref{ux_ricc} and \eqref{HS_eta_a}, we have, whenever
\(D_{\mathrm{HS}}(a,t)\neq0\),
\[
    \partial_a V_{\mathrm{HS}}(a,t)
    =
    v_0'(a)D_{\mathrm{HS}}(a,t)^{\lambda-1}.
\]
Since \(\lambda\geq2\), the right-hand side extends continuously
through \(D_{\mathrm{HS}}=0\).

The differentiated equation \eqref{genHS} is invariant under the
time-dependent Galilean transformation
\[
    v(x,t)\longmapsto v(x-h(t),t)+h'(t).
\]
Thus the identities for \(\zeta_a\) and \(\partial_aV_{\mathrm{HS}}\)
determine the Lagrangian variables only up to a time-dependent spatial
translation. We fix this gauge freedom by imposing
\[
    \zeta(a_0,t)=a_0+t\,v_0(a_0),
    \qquad
    V_{\mathrm{HS}}(a_0,t)=v_0(a_0).
\]
Accordingly, we define the extended Lagrangian map and velocity by
\begin{equation}\label{eta_def}
    \zeta(a,t)
    :=
    a_0+t\,v_0(a_0)
    +
    \int_{a_0}^a
    D_{\mathrm{HS}}(y,t)^\lambda\,dy,
    \qquad t\geq0,
\end{equation}
and
\begin{equation}\label{VHS_def}
    V_{\mathrm{HS}}(a,t)
    :=
    v_0(a_0)
    +
    \int_{a_0}^a
    v_0'(y)D_{\mathrm{HS}}(y,t)^{\lambda-1}\,dy,
    \qquad t\geq0.
\end{equation}
At \(t=0\), these definitions recover the initial data:
\[
    \zeta(a,0)=a,
    \qquad
    V_{\mathrm{HS}}(a,0)=v_0(a).
\]
Moreover, since
\( 
    \partial_t D_{\mathrm{HS}}(a,t)^\lambda
    =
    v_0'(a)D_{\mathrm{HS}}(a,t)^{\lambda-1},
\) 
differentiating \eqref{eta_def} with respect to \(t\) gives
\begin{equation}\label{HS_extended_characteristic}
    \partial_t\zeta(a,t)
    =
    V_{\mathrm{HS}}(a,t),
    \qquad t\geq0.
\end{equation}
Thus the extended variables satisfy the characteristic relation for
all \(t\geq0\). 
Hence the extended variables preserve the characteristic relation
\(
    \partial_t\zeta=V_{\mathrm{HS}}
\)
for all \(t\geq0\). On the classical lifespan, after fixing the above
Galilean gauge, these definitions agree with the classical Lagrangian
map and velocity. 

The following lemma records the basic properties of these extensions.

\begin{lemma}[Lagrangian derivatives and monotonicity]
\label{lem_cross_corrected}
Let \(\lambda\in2\mathbb N\), and let \(\zeta\) and
\(V_{\mathrm{HS}}\) be defined by \eqref{eta_def} and
\eqref{VHS_def}, respectively. Then
\begin{equation}\label{eta_a_simple}
    \partial_a\zeta(a,t)
    =
    D_{\mathrm{HS}}(a,t)^\lambda
    =
    \left(
        1+\frac{t}{\lambda}v_0'(a)
    \right)^\lambda
    \geq0,
\end{equation}
and
\begin{equation}\label{VHS_a_simple}
    \partial_aV_{\mathrm{HS}}(a,t)
    =
    v_0'(a)D_{\mathrm{HS}}(a,t)^{\lambda-1}.
\end{equation}

Moreover, there exists \(\veps_{\mathrm{HS}}>0\) such that, for every
\(t\in[0,t_*+\veps_{\mathrm{HS}})\), the map
\(a\mapsto\zeta(a,t)\) is strictly increasing.
\end{lemma}

\begin{proof}
The identities \eqref{eta_a_simple} and \eqref{VHS_a_simple} follow
directly by differentiating \eqref{eta_def} and \eqref{VHS_def} with
respect to \(a\).
By Assumptions~\ref{ass:HS-main}\textup{(H2)--(H3)}, after decreasing
\(\veps_{\mathrm{HS}}>0\) if necessary,
\(D_{\mathrm{HS}}(\cdot,t)\) does not vanish identically on any
nontrivial interval for
\(
    t\in[0,t_*+\veps_{\mathrm{HS}}).
\) 
Since \(\lambda\) is even,
\(D_{\mathrm{HS}}(\cdot,t)^\lambda\geq0\). Hence, for any
\(a_1<a_2\),
\[
    \zeta(a_2,t)-\zeta(a_1,t)
    =
    \int_{a_1}^{a_2}
    D_{\mathrm{HS}}(y,t)^\lambda\,dy
    >0.
\]
Thus \(a\mapsto\zeta(a,t)\) is strictly increasing.
\end{proof}

Consequently, for every
\(t\in[0,t_*+\veps_{\mathrm{HS}})\),
the map \(\zeta(\cdot,t)\) admits a continuous inverse on its image.
We therefore reconstruct the Eulerian profile by
\begin{equation*}
    v(x,t)
    :=
    V_{\mathrm{HS}}
    \bigl(\zeta^{-1}(x,t),t\bigr).
\end{equation*}
Away from the singular set, this reconstruction agrees with the
classical characteristic solution. At points where
\(D_{\mathrm{HS}}=0\), the Lagrangian formulas
\eqref{eta_def}--\eqref{VHS_def} remain regular and provide the
post-breaking continuation of the Eulerian profile.


\begin{remark}
The preceding continuation is closely parallel to the construction
for the Camassa--Holm equation. Indeed, whenever \(v_0'(a)\neq0\),
\eqref{ux_ricc} gives
\(\textstyle 
    \frac{1}{v_x(\zeta(a,t),t)}
    =
    \frac{D_{\mathrm{HS}}(a,t)}{v_0'(a)}.
\)
Thus \(D_{\mathrm{HS}}\) plays, up to the fixed factor \(v_0'(a)\),
the role of the reciprocal slope along characteristics. In the
Camassa--Holm setting, the corresponding quantity is \(Q\), whose
evolution provides the continuation through wave breaking.
\end{remark}

\begin{remark}[Even versus odd \(\lambda\)]
The restriction to even \(\lambda\) is essential for the post-breaking
Lagrangian continuation. Indeed, if \(\lambda\in2\mathbb N+1\), then
\eqref{eta_a_simple} gives
\(\textstyle
    \partial_a\zeta(a_0,t)\ge0
    \quad\text{for }t\le t_*,
    \partial_a\zeta(a_0,t)<0
    \quad\text{for }t>t_*,
\)
whereas
\(
    \partial_aV_{\mathrm{HS}}(a_0,t)\le0
    \text{ for all }t.
\)
Thus the Lagrangian map loses monotonicity immediately after the first
singular time, and distinct Lagrangian trajectories may cross.

This mechanism can be seen more clearly by considering two labels
\(A<B\) and tracing the corresponding trajectories
\(\zeta(A,t)\) and \(\zeta(B,t)\) in the \((\zeta,v)\)-plane.
Since
\(
    \partial_aV_{\mathrm{HS}}
    =
    v_x\,\partial_a\zeta,
\)
the change of sign of \(v_x\) across the first singular time is realized
differently according to the parity of \(\lambda\).

If \(\lambda\) is even, then
\(
    \partial_a\zeta=D_{\mathrm{HS}}^\lambda\ge0,
\)
so the spatial ordering of the labels is preserved:
\(
    \zeta(A,t)<\zeta(B,t).
\)
The sign change is instead carried by
\(\partial_aV_{\mathrm{HS}}\). Hence an initial ordering of the
Lagrangian values, for instance
\(
    v_0(A)>v_0(B),
\) 
may be reversed after the singular time:
\(
    v(\zeta(A,t),t)<v(\zeta(B,t),t).
\)
Such a reversal is compatible with a single-valued continuous Eulerian
profile.

By contrast, if \(\lambda\) is odd, then the sign change occurs in
\(\partial_a\zeta\). After \(t_*\), one may therefore have
\(
    \zeta(A,t)>\zeta(B,t),
\)
while the ordering of the Lagrangian values remains
\(
    v(\zeta(A,t),t)>v(\zeta(B,t),t).
\)
Hence the Lagrangian parametrization can fold over itself: the same
Eulerian position may correspond to more than one Lagrangian label with
different values of \(V_{\mathrm{HS}}\). In particular, the simple
single-valued Eulerian reconstruction used above is no longer available
after trajectory crossing.

This distinction is illustrated in Figure~\ref{lambda_oddeven}.

\begin{figure}[H]
    \centering
    \includegraphics[width=0.6\linewidth]{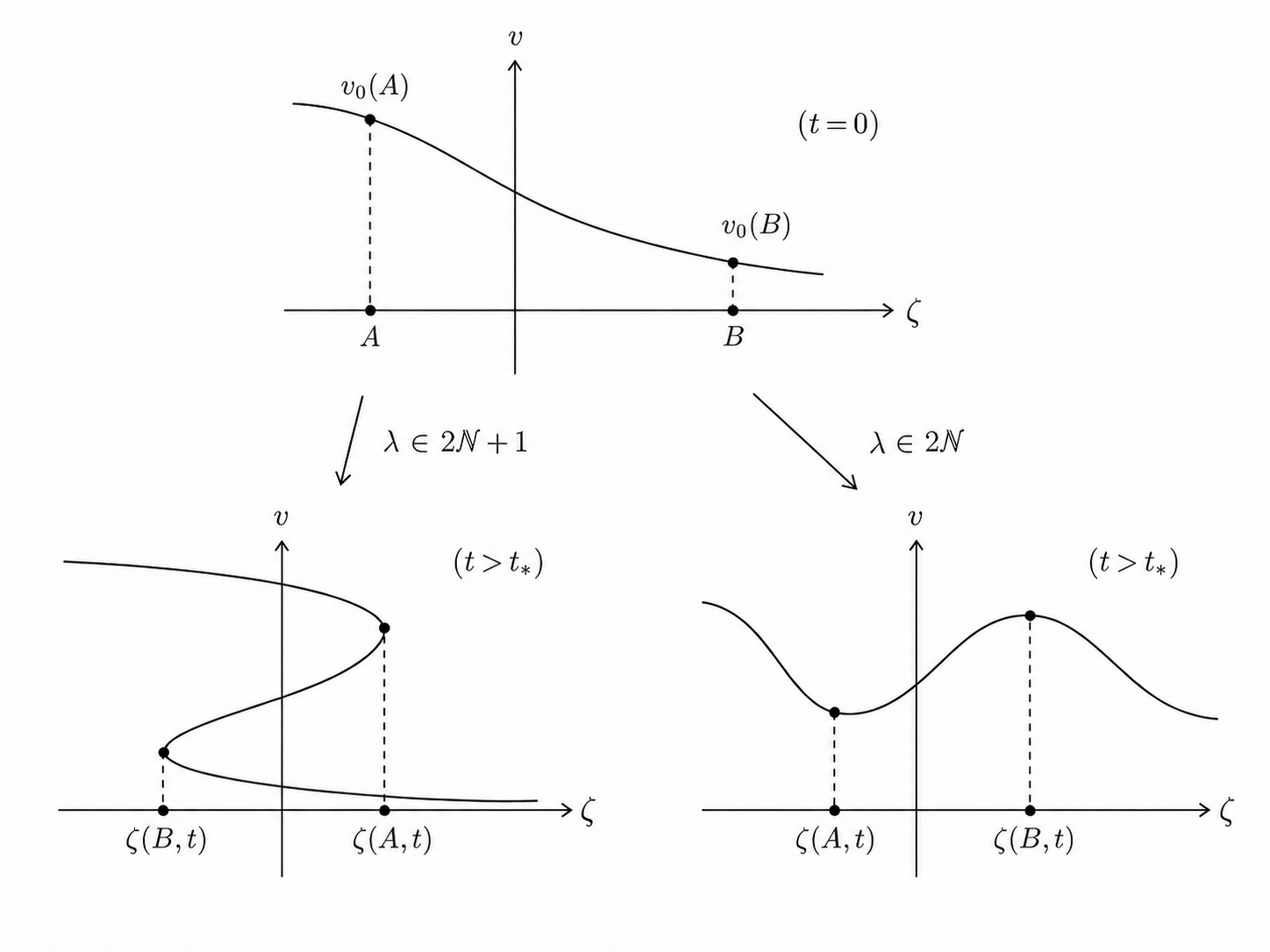}
    \caption{Schematic post-breaking dynamics for even and odd
    values of \(\lambda\).}
    \label{lambda_oddeven}
\end{figure}

For this reason, throughout the paper we restrict attention to
\(\lambda\in2\mathbb N\).
\end{remark}

We first describe the local geometry of the Lagrangian singular set
\(\mathcal S_L\) near the first breaking point.

\begin{lemma}[Local structure of the Lagrangian singular set]
\label{lem:singular_set}
Near \(a=a_0\), the Lagrangian singular set
\(\mathcal S_L\) defined in \eqref{lagrangian_singular_set} is given by
the smooth graph
\[
    t=t_s(a),
    \qquad
    t_s(a):=-\frac{\lambda}{v_0'(a)}.
\]
At the first singular point \((a_0,t_*)\),
\(
    t_s'(a_0)=0,
\)
and
\begin{equation}\label{ts''}
    t_s''(a_0)
    =
    \lambda
    \frac{v_0^{(3)}(a_0)}
         {v_0'(a_0)^2}
    >0.
\end{equation}
Hence \(t_s\) has a nondegenerate local minimum at \(a=a_0\).

Moreover, the Eulerian projection of the singular curve,
\(
    X_s(a):=\zeta(a,t_s(a)),
\)
is well defined near \(a=a_0\), and
\(
    X_s'(a_0)=0.
\)
\end{lemma}

\begin{proof}
Since
\(
    v_0'(a_0)<0,
\)
the function
\(\textstyle 
    t_s(a)
    =
    -\frac{\lambda}{v_0'(a)}
\)
is smooth in a neighborhood of \(a=a_0\). By
\(D_{\mathrm{HS}}(a,t)=0\) if and only if \(t=t_s(a)\),
this gives a local graph representation of \(\mathcal S_L\).

Differentiating \(t_s\), we obtain
\[
    t_s'(a)
    =
    \lambda
    \frac{v_0''(a)}
         {v_0'(a)^2}.
\]
Since
\(
    v_0''(a_0)=0,
\)
it follows that
\(
    t_s'(a_0)=0.
\)
Differentiating once more gives
\[
    t_s''(a)
    =
    \lambda
    \left(
        \frac{v_0^{(3)}(a)}
             {v_0'(a)^2}
        -
        2
        \frac{v_0''(a)^2}
             {v_0'(a)^3}
    \right).
\]
Evaluating at \(a=a_0\) and using
\(
    v_0''(a_0)=0,
    v_0^{(3)}(a_0)>0,
\) 
we obtain \eqref{ts''}. Thus \(t_s\) has a nondegenerate local
minimum at \(a=a_0\). In particular, the Lagrangian singular curve
has a horizontal tangent at the first breaking point
\((a_0,t_*)\).

We next consider its Eulerian projection
\(
    X_s(a)=\zeta(a,t_s(a)).
\)
By the chain rule,
\(
    X_s'(a)
    =
    \zeta_a(a,t_s(a))
    +
    \zeta_t(a,t_s(a))\,t_s'(a).
\)
At \(a=a_0\), we have
\(
    \zeta_a(a_0,t_*)=0
\)
by Lemma~\ref{lem_cross_corrected}, while
\(
    t_s'(a_0)=0.
\)
Therefore,
\(
    X_s'(a_0)=0.
\)
This completes the proof.
\end{proof}
The nondegenerate minimum of \(t_s\) immediately yields the
bifurcation of the singular label \(a_0\) into two branches
\(a_-(t)\) and \(a_+(t)\) for \(t>t_*\) sufficiently close to
\(t_*\).
\subsection{Emergence of a pair of cuspons}
\label{sec_HScusp}

We now describe the short-time post-breaking geometry and show that a pair
of cusps emerges from the first breaking point.

A useful feature of the generalized Hunter--Saxton equation is that its
Lagrangian degeneracy is completely explicit. Indeed, by
\eqref{eta_a_simple},
\[
    \zeta_a(a,t)=D_{\mathrm{HS}}(a,t)^\lambda,
    \qquad
    D_{\mathrm{HS}}(a,t)
    =
    1+\frac{t}{\lambda}v_0'(a).
\]
Thus the singular labels are characterized simply by
\begin{equation}\label{HS_singular_label_equation}
    D_{\mathrm{HS}}(a,t)=0,
    \qquad\text{equivalently}\qquad
    v_0'(a)=-\frac{\lambda}{t}.
\end{equation}
This explicit characterization is in contrast with the Camassa--Holm
case, where the corresponding singular labels are governed by the
nonlocal Lagrangian dynamics.

By \eqref{HS_tstar}, the first singularity occurs at the label-time pair \((a_0, t_*) \), where 
\(\textstyle 
    t_*=-\frac{\lambda}{v_0'(a_0)}.
\) 
Moreover, after decreasing \(\veps_{\mathrm{HS}}>0\) if necessary, for
every
\(
    t\in(t_*,t_*+\veps_{\mathrm{HS}})
\)
the level equation \eqref{HS_singular_label_equation} has exactly two
solutions near \(a_0\), which we denote by
\begin{equation*}
    a_-(t)<a_0<a_+(t).
\end{equation*}
Equivalently,
\begin{equation*}
    \zeta_a(a_\pm(t),t)=0.
\end{equation*}
The two branches \(a_\pm(t)\) bifurcate from the first breaking label:
\[
    a_\pm(t)\longrightarrow a_0
    \qquad\text{as }t\downarrow t_*.
\]

Since no confusion with the Camassa--Holm variables will arise, we use
the same notation \(a_\pm(t)\) for the two post-breaking singular labels
and define their Eulerian images by
\begin{equation*}
    x_\pm(t)
    :=
    \zeta(a_\pm(t),t).
\end{equation*}
Thus we also use the same notation \(x_\pm(t)\) as in the
Camassa--Holm setting. We shall show below that \(x_-(t)\) and
\(x_+(t)\) are precisely the locations of the two cusps. 
\begin{lemma}[Reversal of the Lagrangian profile and emergence of a cusp pair]
\label{HS_blowup}
There exists \(\veps_{\mathrm{HS}}>0\) such that the following hold.

For \(t<t_*\), the Lagrangian profile
\(
    V_{\mathrm{HS}}(a,t)
    :=
    v(\zeta(a,t),t)
\)
is strictly decreasing in a neighborhood of \(a=a_0\). In particular,
for \(a,b\) sufficiently close to \(a_0\) with
\(
    b<a_0<a,
\)
we have
\begin{equation}\label{pre_breaking_ordering}
    v(\zeta(b,t),t)
    >
    v(\zeta(a_0,t),t)
    >
    v(\zeta(a,t),t).
\end{equation}

By contrast, for every
\(
    t\in(t_*,t_*+\veps_{\mathrm{HS}}),
\)
the Lagrangian profile is strictly increasing on
\(
    (a_-(t),a_+(t)).
\)
Consequently, for
\(
    a_-(t)<b<a_0<a<a_+(t),
\)
we have
\begin{equation}\label{post_breaking_ordering}
    v(\zeta(b,t),t)
    <
    v(\zeta(a_0,t),t)
    <
    v(\zeta(a,t),t).
\end{equation}

Moreover,
\[
    \left|
        v_x(\zeta(a,t),t)
    \right|
    \longrightarrow\infty
    \qquad
    \text{as }a\to a_\pm(t).
\]
Hence the continuous Eulerian profile develops two singular points at
\(
    x_\pm(t)=\zeta(a_\pm(t),t),
\)
with unbounded spatial slope. 
\end{lemma}

\begin{proof}
By Lemma~\ref{lem_cross_corrected},
\begin{equation}\label{HS_Va_sign_formula}
    \partial_a V_{\mathrm{HS}}(a,t)
    =
    v_0'(a)D_{\mathrm{HS}}(a,t)^{\lambda-1},
\end{equation}
where
\(\textstyle 
    D_{\mathrm{HS}}(a,t)
    =
    1+\frac{t}{\lambda}v_0'(a).
\) 

For every \(t<t_*\), the definition of the first singular time gives
\(
    D_{\mathrm{HS}}(a,t)>0
    \text{ for all }a\in\mathbb R.
\) 
Since
\(
    v_0'(a_0)<0,
\)
continuity of \(v_0'\) yields a neighborhood of \(a_0\) on which
\(v_0'<0\). Because \(\lambda-1\) is odd, it follows from
\eqref{HS_Va_sign_formula} that
\(
    \partial_aV_{\mathrm{HS}}(a,t)<0
\)
in this neighborhood. Thus \(V_{\mathrm{HS}}(\cdot,t)\) is strictly
decreasing near \(a_0\), which proves
\eqref{pre_breaking_ordering}.

We now fix
\(
    t\in(t_*,t_*+\veps_{\mathrm{HS}}).
\)
By the definition of the two post-breaking singular labels,
\(
    D_{\mathrm{HS}}(a_\pm(t),t)=0,
\)
or equivalently,
\(\textstyle 
    v_0'(a_\pm(t))
    =
    -\frac{\lambda}{t}.
\)
On the other hand,
\(\textstyle 
    v_0'(a_0)
    =
    -\frac{\lambda}{t_*}
    <
    -\frac{\lambda}{t},
\) 
since \(t>t_*\). Hence
\(
    D_{\mathrm{HS}}(a_0,t)<0.
\)
For \(\veps_{\mathrm{HS}}>0\) sufficiently small,
\(a_-(t)\) and \(a_+(t)\) are the only zeros of
\(D_{\mathrm{HS}}(\cdot,t)\) in
\([a_-(t),a_+(t)]\). Therefore,
\(
    D_{\mathrm{HS}}(a,t)<0
    \text{ for }
    a\in(a_-(t),a_+(t)).
\) 
In particular,
\(\textstyle 
    v_0'(a)
    <
    -\frac{\lambda}{t}
    <
    0
    \text{ for }
    a\in(a_-(t),a_+(t)).
\) 
Since \(\lambda-1\) is odd, both factors on the right-hand side of
\eqref{HS_Va_sign_formula} are negative, and hence
\[
    \partial_aV_{\mathrm{HS}}(a,t)>0
    \qquad
    \text{for }
    a\in(a_-(t),a_+(t)).
\]
Thus \(V_{\mathrm{HS}}(\cdot,t)\) is strictly increasing on this
interval, and \eqref{post_breaking_ordering} follows.

Finally, by \eqref{ux_ricc},
\begin{equation*}
    v_x(\zeta(a,t),t)
    =
    \frac{v_0'(a)}
         {D_{\mathrm{HS}}(a,t)}
\end{equation*}
whenever \(D_{\mathrm{HS}}(a,t)\neq0\).
Since
\( \textstyle 
    D_{\mathrm{HS}}(a_\pm(t),t)=0,
    v_0'(a_\pm(t))
    =
    -\frac{\lambda}{t}\neq0,
\)
we obtain
\[
    \left|
        v_x(\zeta(a,t),t)
    \right|
    \longrightarrow\infty
    \qquad
    \text{as }a\to a_\pm(t).
\]
Since \(a\mapsto\zeta(a,t)\) remains strictly increasing and the
Eulerian profile reconstructed above is continuous, the points
\(
    x_\pm(t)=\zeta(a_\pm(t),t)
\)
are two distinct singular points of the continuous profile with
unbounded spatial slope. The precise cusp asymptotics at these points
will be established below.
\end{proof}
By \eqref{ux_ricc} and the uniqueness of the global minimum of \(v_0'\),
the slope becomes singular at the first singular time \(t=t_*\) only at
the first breaking point
\(
    x_*:=\zeta(a_0,t_*).
\)
Lemma~\ref{HS_blowup} shows that, immediately after \(t_*\), the
corresponding singular label bifurcates into the two branches
\(a_-(t)\) and \(a_+(t)\), whose Eulerian images are
\(
    x_\pm(t):=\zeta(a_\pm(t),t).
\)
Thus the single singular point \(x_*\) gives rise to a pair of
post-breaking singular points.

The sign structure also determines their orientations. For
\(t\in(t_*,t_*+\veps_{\mathrm{HS}})\), we have
\[
    v_x(\zeta(a,t),t)>0
    \qquad
    \text{for }a\in(a_-(t),a_+(t)),
\]
whereas \(v_x(\zeta(a,t),t)<0\) immediately outside this interval.
Consequently, \(x_-(t)\) is of local-minimum type and \(x_+(t)\) is of
local-maximum type. This reversal of the local profile is the basic
mechanism underlying the formation of the cusp pair: before breaking,
the Lagrangian profile is decreasing near \(a=a_0\), while immediately
after breaking it becomes increasing between the two singular labels.

We next record that the continued solution remains classical away from
the two singular curves.

\begin{lemma}[Classical regularity away from the cusp pair]
\label{lem_class}
After decreasing \(\veps_{\mathrm{HS}}>0\) if necessary, define
\[
    \Omega_{\mathrm{HS}}
    :=
    \left\{
        (x,t)\in\mathbb R\times
        (t_*,t_*+\veps_{\mathrm{HS}})
        :
        x\neq x_-(t),\;
        x\neq x_+(t)
    \right\}.
\]
Then, for every
\(t\in(t_*,t_*+\veps_{\mathrm{HS}})\), it holds that 
\(
    v(\cdot,t)
    \in
    C^1\bigl(
        \mathbb R\setminus\{x_-(t),x_+(t)\}
    \bigr).
\)
Moreover,
\(
    v_x\in C^1(\Omega_{\mathrm{HS}}),
\)
and \(v\) satisfies the generalized Hunter--Saxton equation
\eqref{genHS} pointwise on \(\Omega_{\mathrm{HS}}\).
\end{lemma}

\begin{proof}
Let
\(
    (x,t)\in\Omega_{\mathrm{HS}}.
\)
Since \(a\mapsto\zeta(a,t)\) is strictly increasing, there exists a
unique label \(a\neq a_\pm(t)\) such that
\(
    x=\zeta(a,t).
\)
In particular,
\(
    D_{\mathrm{HS}}(a,t)\neq0.
\)
By Lemma~\ref{lem_cross_corrected},
\(
    \zeta_a(a,t)
    =
    D_{\mathrm{HS}}(a,t)^\lambda
    \neq0.
\)
Hence the map
\(
    (a,t)\longmapsto(\zeta(a,t),t)
\)
is locally invertible on the preimage of
\(\Omega_{\mathrm{HS}}\). Since \(D_{\mathrm{HS}}\neq0\) there, the
extended Lagrangian variables \(\zeta\) and \(V_{\mathrm{HS}}\) are
smooth on this region. The inverse function theorem therefore gives
a smooth inverse Lagrangian label
\(
    a=A(x,t)
\) on \(\Omega_{\mathrm{HS}}\), where $A(x,t):=(\zeta(\cdot,t))^{-1}(x)$.

By the Eulerian reconstruction,
\(
    v(x,t)
    =
    V_{\mathrm{HS}}(A(x,t),t),
\)
and hence
\(
    v\in C^1(\Omega_{\mathrm{HS}}).
\) 
Differentiating with respect to \(x\) and using
\eqref{eta_a_simple}--\eqref{VHS_a_simple}, we obtain
\begin{equation}\label{HS_vx_eulerian}
    v_x(x,t)
    =
    \frac{
        \partial_aV_{\mathrm{HS}}(A(x,t),t)
    }{
        \partial_a\zeta(A(x,t),t)
    }
    =
    \frac{
        v_0'(A(x,t))
    }{
        D_{\mathrm{HS}}(A(x,t),t)
    }.
\end{equation}
Since \(D_{\mathrm{HS}}\neq0\) on the region under consideration and
\(A\) is smooth, the right-hand side of
\eqref{HS_vx_eulerian} shows that
\(
    v_x\in C^1(\Omega_{\mathrm{HS}}).
\) 

For later use, we also record the second derivative explicitly.
Differentiating
\(\textstyle 
    \frac{v_0'(a)}{D_{\mathrm{HS}}(a,t)}
\)
with respect to \(a\), and recalling that
\(\textstyle 
    D_{\mathrm{HS}}(a,t)
    =
    1+\frac{t}{\lambda}v_0'(a),
\)
gives
\(\textstyle 
    \partial_a
    \left(
        \frac{v_0'(a)}
             {D_{\mathrm{HS}}(a,t)}
    \right)
    =
    \frac{v_0''(a)}
         {D_{\mathrm{HS}}(a,t)^2}.
\)
Since
\(
    \zeta_a(a,t)
    =
    D_{\mathrm{HS}}(a,t)^\lambda,
\)
we therefore obtain
\begin{equation*}
    v_{xx}(\zeta(a,t),t)
    =
    \frac{v_0''(a)}
         {D_{\mathrm{HS}}(a,t)^{\lambda+2}},
\end{equation*}
whenever \(a\neq a_\pm(t)\). Thus \(v_{xx}\) is continuous on
\(\Omega_{\mathrm{HS}}\).

It remains to verify the equation. By
\eqref{eta_def}--\eqref{VHS_def},
\[
    \partial_t\zeta(a,t)
    =
    V_{\mathrm{HS}}(a,t)
    =
    v(\zeta(a,t),t),
\]
so the extended Lagrangian map continues to satisfy the characteristic
equation. Moreover, \eqref{HS_vx_eulerian} gives
\(\textstyle 
    v_x(\zeta(a,t),t)
    =
    \frac{v_0'(a)}
         {D_{\mathrm{HS}}(a,t)}.
\)
Differentiating with respect to \(t\) at fixed \(a\), we find
\[
\begin{aligned}
    \frac{d}{dt}
    v_x(\zeta(a,t),t)
    &=
    -\frac{1}{\lambda}
    \frac{v_0'(a)^2}
         {D_{\mathrm{HS}}(a,t)^2}  =
    -\frac{1}{\lambda}
    v_x(\zeta(a,t),t)^2.
\end{aligned}
\]
On the other hand, since
\(
    \partial_t\zeta(a,t)
    =
    v(\zeta(a,t),t),
\) 
the chain rule yields
\[
    \frac{d}{dt}
    v_x(\zeta(a,t),t)
    =
    v_{tx}(\zeta(a,t),t)
    +
    v(\zeta(a,t),t)
    v_{xx}(\zeta(a,t),t).
\]
Combining the last two identities, we obtain
\[
    v_{tx}
    +
    v v_{xx}
    =
    -\frac{1}{\lambda}v_x^2
\]
at every point of \(\Omega_{\mathrm{HS}}\). Hence \(v\) is a
classical solution of the generalized Hunter--Saxton equation
\eqref{genHS} on \(\Omega_{\mathrm{HS}}\).
\end{proof}

%
\subsection{Sharp regularity and dynamics of the cuspon pair}
\label{sec_HS_reg}

We first determine the sharp spatial regularity at the first breaking
point and at the two post-breaking cusp locations.

\begin{lemma}[Sharp regularity of the cusp pair]
\label{reg_lem}
There exists \(\veps_{\mathrm{HS}}>0\) such that the following hold.

At the first singular time \(t=t_*\), the pointwise H\"older exponent
of \(v(\cdot,t_*)\) at
\(
    x_*=\zeta(a_0,t_*)
\)
is exactly
\begin{equation*}
    \alpha_*
    :=
    \frac{2\lambda-1}{2\lambda+1}.
\end{equation*}
More precisely,
\begin{equation}\label{HS_holder_first_limit}
    0
    <
    \lim_{x\to x_*}
    \frac{
        |v(x,t_*)-v(x_*,t_*)|
    }{
        |x-x_*|^{\alpha_*}
    }
    <
    \infty.
\end{equation}

For every
\(
    t\in(t_*,t_*+\veps_{\mathrm{HS}}),
\)
the pointwise H\"older exponent of \(v(\cdot,t)\) at each of the two
post-breaking singular points
\(
    x_\pm(t)=\zeta(a_\pm(t),t)
\)
is exactly
\begin{equation*}
    \alpha_{\mathrm{cusp}}
    :=
    \frac{\lambda}{\lambda+1}.
\end{equation*}
More precisely, for either choice of sign,
\begin{equation}\label{HS_holder_post_limit}
    0
    <
    \lim_{x\to x_\pm(t)}
    \frac{
        |v(x,t)-v(x_\pm(t),t)|
    }{
        |x-x_\pm(t)|^{\alpha_{\mathrm{cusp}}}
    }
    <
    \infty.
\end{equation}
In particular, both exponents are sharp.
\end{lemma}

\begin{proof}
We treat the first breaking point and the post-breaking cusp points
simultaneously.

Fix \(t\in[t_*,t_*+\veps_{\mathrm{HS}})\), and let \(a_s\) denote one
of the singular labels under consideration:
\[
    a_s
    :=
    \begin{cases}
        a_0,
        & t=t_*,
        \\[1mm]
        a_-(t)\ \text{or}\ a_+(t),
        & t>t_*.
    \end{cases}
\]
Set
\(
    x_s:=\zeta(a_s,t).
\)
We also define
\[
    m
    :=
    \begin{cases}
        2,
        & t=t_*,
        \\[1mm]
        1,
        & t>t_*.
    \end{cases}
\]

By the definition of a singular label,
\begin{equation}\label{HS_singular_slope_relation}
    D_{\mathrm{HS}}(a_s,t)=0,
    \qquad
    v_0'(a_s)=-\frac{\lambda}{t}.
\end{equation}
At \(t=t_*\), the nondegeneracy of the minimum of \(v_0'\) gives
\(
    v_0''(a_0)=0,
    v_0^{(3)}(a_0)\neq0,
\)
whereas for \(t>t_*\) sufficiently close to \(t_*\), the two level
crossings \(a_\pm(t)\) are simple, so that
\(
    v_0''(a_\pm(t))\neq0.
\)
Consequently, as \(\alpha\to a_s\),
\begin{equation}\label{HS_D_local_expansion}
    D_{\mathrm{HS}}(\alpha,t)
    =
    d_s(t)(\alpha-a_s)^m
    +
    O\left(|\alpha-a_s|^{m+1}\right),
\end{equation}
where
\[
    d_s(t)
    :=
    \frac{t}{\lambda}
    \frac{v_0^{(m+1)}(a_s)}{m!}
    \neq0.
\]

By \eqref{eta_a_simple},
\(
    \zeta_a(a,t)
    =
    D_{\mathrm{HS}}(a,t)^\lambda.
\)
Hence, using \eqref{HS_D_local_expansion},
\begin{align}
    \zeta(a,t)-\zeta(a_s,t)
    &=
    \int_{a_s}^a
    D_{\mathrm{HS}}(\alpha,t)^\lambda\,d\alpha
    \notag
    =
    \frac{d_s(t)^\lambda}{m\lambda+1}
    (a-a_s)^{m\lambda+1}
    +
    O\left(
        |a-a_s|^{m\lambda+2}
    \right).
    \notag
\end{align}
Setting
\(\textstyle 
    C_s(t)
    :=
    \frac{d_s(t)^\lambda}{m\lambda+1},
\) 
we obtain
\begin{equation}\label{HS_zeta_cusp_expansion_short}
    \zeta(a,t)-x_s
    =
    C_s(t)(a-a_s)^{m\lambda+1}
    +
    O\left(
        |a-a_s|^{m\lambda+2}
    \right),
\end{equation}
where \(C_s(t)\neq0\).

We next expand the Lagrangian profile. By
\eqref{VHS_a_simple}, we see that
\(
    \partial_aV_{\mathrm{HS}}(a,t)
    =
    v_0'(a)D_{\mathrm{HS}}(a,t)^{\lambda-1}.
\)
Using \eqref{HS_singular_slope_relation} and
\eqref{HS_D_local_expansion}, we obtain
\begin{align}
    V_{\mathrm{HS}}(a,t)-V_{\mathrm{HS}}(a_s,t)
    &=
    \int_{a_s}^a
    v_0'(\alpha)
    D_{\mathrm{HS}}(\alpha,t)^{\lambda-1}
    \,d\alpha
    \notag\\
    &=
    -\frac{\lambda}{t}
    \frac{
        d_s(t)^{\lambda-1}
    }{
        m(\lambda-1)+1
    }
    (a-a_s)^{m(\lambda-1)+1}
    \notag
    +
    O\left(
        |a-a_s|^{m(\lambda-1)+2}
    \right).
    \notag
\end{align}
Let
\(\textstyle 
    \widetilde C_s(t)
    :=
    -\frac{\lambda}{t}
    \frac{
        d_s(t)^{\lambda-1}
    }{
        m(\lambda-1)+1
    }.
\) 
Then \(\widetilde C_s(t)\neq0\), and
\begin{equation}\label{HS_V_cusp_expansion_short}
    V_{\mathrm{HS}}(a,t)-V_{\mathrm{HS}}(a_s,t)
    =
    \widetilde C_s(t)
    (a-a_s)^{m(\lambda-1)+1}
    +
    O\left(
        |a-a_s|^{m(\lambda-1)+2}
    \right).
\end{equation}

Combining \eqref{HS_zeta_cusp_expansion_short} and
\eqref{HS_V_cusp_expansion_short}, we find that the critical
H\"older exponent at \(x_s\) is
\[
    \alpha_s
    =
    \frac{
        m(\lambda-1)+1
    }{
        m\lambda+1
    }.
\]
Indeed,
\begin{equation}\label{HS_holder_ratio_limit_general}
    \frac{
        |V_{\mathrm{HS}}(a,t)-V_{\mathrm{HS}}(a_s,t)|
    }{
        |\zeta(a,t)-\zeta(a_s,t)|^{\alpha_s}
    }
    \longrightarrow
    \frac{
        |\widetilde C_s(t)|
    }{
        |C_s(t)|^{\alpha_s}
    }
    \in(0,\infty)
\end{equation}
as \(a\to a_s\).

At the first singular time \(t=t_*\), we have \(m=2\), and therefore
\[
    \alpha_s
    =
    \frac{2(\lambda-1)+1}{2\lambda+1}
    =
    \frac{2\lambda-1}{2\lambda+1}.
\]
This proves \eqref{HS_holder_first_limit}.

For every \(t>t_*\) sufficiently close to \(t_*\), each singular
label \(a_\pm(t)\) is a simple zero of
\(D_{\mathrm{HS}}(\cdot,t)\), so that \(m=1\). Hence
\[
    \alpha_s
    =
    \frac{\lambda}{\lambda+1},
\]
which proves \eqref{HS_holder_post_limit}
at both \(x_-(t)\) and \(x_+(t)\).

The positivity and finiteness of the limit in
\eqref{HS_holder_ratio_limit_general} also show that these exponents
cannot be improved. Thus the stated H\"older exponents are sharp.
\end{proof}
\begin{remark}[Improved regularity after cusp bifurcation]
An interesting feature of the post-breaking dynamics is that the two
cuspons are more regular than the profile at the first breaking point.
Indeed,
\(\textstyle 
    \frac{\lambda}{\lambda+1}
    >
    \frac{2\lambda-1}{2\lambda+1}.
\)
Thus, at \(t=t_*\), the first singular profile has the sharp pointwise
H\"older exponent
\(\textstyle 
    \frac{2\lambda-1}{2\lambda+1},
\)
whereas immediately after breaking each of the two cuspons has the
strictly larger exponent
\(\textstyle 
    \frac{\lambda}{\lambda+1}.
\)
In this sense, the bifurcation of the single breaking singularity into
an anticuspon--cuspon pair is accompanied by an improvement of the
pointwise spatial regularity at the singular points.
\end{remark}
\begin{lemma}[Outward motion of the cusp pair]
\label{a+far_lem}
After decreasing \(\veps_{\mathrm{HS}}>0\) if necessary, let
\(
    t\in(t_*,t_*+\veps_{\mathrm{HS}}).
\)
Then the two cusp locations
\(
    x_\pm(t):=\zeta(a_\pm(t),t)
\)
move away from the central trajectory
\(
    x_0(t):=\zeta(a_0,t).
\)
More precisely,
\[
    \frac{d}{dt}\bigl(x_+(t)-x_0(t)\bigr)>0,
    \qquad
    \frac{d}{dt}\bigl(x_0(t)-x_-(t)\bigr)>0.
\]
Consequently,
\[
    \frac{d}{dt}\bigl(x_+(t)-x_-(t)\bigr)>0.
\]
In fact,
\[
    x_-'(t)<x_0'(t)<x_+'(t).
\]
\end{lemma}

\begin{proof}
Since \(a_\pm(t)\) are simple zeros of
\(D_{\mathrm{HS}}(\cdot,t)\), they depend smoothly on \(t\) for
\(t>t_*\) sufficiently close to \(t_*\). Moreover,
\(
    \zeta_a(a_\pm(t),t)=0.
\)
By \eqref{HS_extended_characteristic},
\[
\begin{aligned}
    \frac{d}{dt}x_\pm(t)
    &=
    \zeta_t(a_\pm(t),t)
    +
    \zeta_a(a_\pm(t),t)a_\pm'(t) =
    V_{\mathrm{HS}}(a_\pm(t),t).
\end{aligned}
\] 
Since \(a_0\) is fixed,
\[
    x_0'(t)
    =
    \zeta_t(a_0,t)
    =
    V_{\mathrm{HS}}(a_0,t).
\]
Hence
\[
    \frac{d}{dt}\bigl(x_+(t)-x_0(t)\bigr)
    =
    V_{\mathrm{HS}}(a_+(t),t)
    -
    V_{\mathrm{HS}}(a_0,t).
\]

By the sign property established in the proof of
Lemma~\ref{HS_blowup},
\[
    \partial_aV_{\mathrm{HS}}(a,t)>0
    \qquad
    \text{for }
    a\in(a_-(t),a_+(t)).
\]
Thus
\[
    V_{\mathrm{HS}}(a_+(t),t)
    -
    V_{\mathrm{HS}}(a_0,t)
    =
    \int_{a_0}^{a_+(t)}
    \partial_aV_{\mathrm{HS}}(a,t)\,da
    >0,
\]
and therefore
\[
    \frac{d}{dt}\bigl(x_+(t)-x_0(t)\bigr)>0.
\]

Similarly,
\[
    \frac{d}{dt}\bigl(x_0(t)-x_-(t)\bigr)
    =
    V_{\mathrm{HS}}(a_0,t)
    -
    V_{\mathrm{HS}}(a_-(t),t)
\]
and hence
\[
    V_{\mathrm{HS}}(a_0,t)
    -
    V_{\mathrm{HS}}(a_-(t),t)
    =
    \int_{a_-(t)}^{a_0}
    \partial_aV_{\mathrm{HS}}(a,t)\,da
    >0.
\]
Therefore,
\[
    \frac{d}{dt}\bigl(x_0(t)-x_-(t)\bigr)>0.
\]

Finally,
\(
    x_+(t)-x_-(t)
    =
    \bigl(x_+(t)-x_0(t)\bigr)
    +
    \bigl(x_0(t)-x_-(t)\bigr),
\)
so
\[
    \frac{d}{dt}\bigl(x_+(t)-x_-(t)\bigr)>0.
\]
This completes the proof.
\end{proof}

\subsection{Completion of the proof of Theorem~\ref{mainthm1}}
\label{sec_thmpf}

\begin{proof}[Proof of \eqref{mainthm_space'} and \eqref{mainthm_space}
in Theorem~\ref{mainthm1}]

Since \(\lambda\) is even, Lemma~\ref{lem_cross_corrected} gives
\(
    \zeta_a(a,t)
    =
    D_{\mathrm{HS}}(a,t)^\lambda
    \geq0.
\)
Moreover, for
\(
    t\in[0,t_*+\veps_{\mathrm{HS}}),
\)
the function \(D_{\mathrm{HS}}(\cdot,t)\) does not vanish identically
on any nontrivial interval. Hence, for any \(a_1<a_2\),
\[
    \zeta(a_2,t)-\zeta(a_1,t)
    =
    \int_{a_1}^{a_2}
    D_{\mathrm{HS}}(a,t)^\lambda\,da
    >0.
\]
Thus \(a\mapsto\zeta(a,t)\) is strictly increasing and admits a
continuous inverse on its image. Consequently, the Eulerian profile
\(
    v(x,t)
    =
    V_{\mathrm{HS}}\bigl(\zeta^{-1}(x,t),t\bigr)
\)
is well defined and continuous for
\(t\in[0,t_*+\veps_{\mathrm{HS}})\).

We next prove conservation of the \(L^\lambda\)-norm of the slope.
By \eqref{eta_a_simple} and \eqref{VHS_a_simple},
\[
    \zeta_a(a,t)
    =
    D_{\mathrm{HS}}(a,t)^\lambda,
    \qquad
    \partial_aV_{\mathrm{HS}}(a,t)
    =
    v_0'(a)D_{\mathrm{HS}}(a,t)^{\lambda-1}.
\]
Away from the singular labels,
\[
    v_x(\zeta(a,t),t)
    =
    \frac{\partial_aV_{\mathrm{HS}}(a,t)}
         {\zeta_a(a,t)}
    =
    \frac{v_0'(a)}
         {D_{\mathrm{HS}}(a,t)}.
\]
Therefore, using the change of variables
\(x=\zeta(a,t)\),
\[
    \begin{aligned}
    \int_{\mathbb R}|v_x(x,t)|^\lambda\,dx
    &=
    \int_{\mathbb R}
    \left|
        \frac{\partial_aV_{\mathrm{HS}}(a,t)}
             {\zeta_a(a,t)}
    \right|^\lambda
    \zeta_a(a,t)\,da
    =
    \int_{\mathbb R}
    \frac{
        |\partial_aV_{\mathrm{HS}}(a,t)|^\lambda
    }{
        \zeta_a(a,t)^{\lambda-1}
    }
    \,da
    =
    \int_{\mathbb R}
    |v_0'(a)|^\lambda\,da.
    \end{aligned}
\]
Hence
\(
    \|v_x(\cdot,t)\|_{L^\lambda(\mathbb R)}
    =
    \|v_0'\|_{L^\lambda(\mathbb R)}
\)
for every
\(t\in[0,t_*+\veps_{\mathrm{HS}})\).

We also obtain an \(L^1\)-bound for the slope. By
\eqref{VHS_a_simple},
\[
    |\partial_aV_{\mathrm{HS}}(a,t)|
    =
    |v_0'(a)|
    \left|
        1+\frac{t}{\lambda}v_0'(a)
    \right|^{\lambda-1}.
\]
Using
\(
    |1+z|^{\lambda-1}
    \leq
    2^{\lambda-1}
    \bigl(1+|z|^{\lambda-1}\bigr),
\)
we obtain
\[
    |\partial_aV_{\mathrm{HS}}(a,t)|
    \leq
    2^{\lambda-1}
    \left(
        |v_0'(a)|
        +
        \left(\frac{t}{\lambda}\right)^{\lambda-1}
        |v_0'(a)|^\lambda
    \right).
\]
Since
\(
    v_0'\in L^1(\mathbb R)\cap L^\lambda(\mathbb R),
\)
it follows that
\(
    \partial_aV_{\mathrm{HS}}(\cdot,t)
    \in L^1(\mathbb R).
\)
Again using the monotone change of variables,
\[
    \begin{aligned}
    \int_{\mathbb R}|v_x(x,t)|\,dx
    &=
    \int_{\mathbb R}
    \left|
        \frac{\partial_aV_{\mathrm{HS}}(a,t)}
             {\zeta_a(a,t)}
    \right|
    \zeta_a(a,t)\,da
    =
    \int_{\mathbb R}
    |\partial_aV_{\mathrm{HS}}(a,t)|\,da
    <\infty.
    \end{aligned}
\]
Thus
\(
    v_x(\cdot,t)\in L^1(\mathbb R).
\)

Finally, we prove the local boundedness of \(v\).
For any \(a,b\in\mathbb R\),
\[
    V_{\mathrm{HS}}(a,t)-V_{\mathrm{HS}}(b,t)
    =
    \int_b^a
    \partial_\alpha V_{\mathrm{HS}}(\alpha,t)\,d\alpha.
\]
Therefore,
\[
    |V_{\mathrm{HS}}(a,t)|
    \leq
    |V_{\mathrm{HS}}(a_0,t)|
    +
    \|\partial_aV_{\mathrm{HS}}(\cdot,t)\|_{L^1(\mathbb R)}.
\]
By the construction of the extended Lagrangian flow,
\(V_{\mathrm{HS}}(a_0,t)=\zeta_t(a_0,t)\) is continuous in \(t\).
Hence it remains bounded on every compact time interval.
Moreover, the preceding \(L^1\)-estimate is uniform for
\(t\) in compact subsets of
\([0,t_*+\veps_{\mathrm{HS}})\).
Consequently,
\(
    v\in
    L^\infty_{\mathrm{loc}}
    \bigl(
        \mathbb R\times[0,t_*+\veps_{\mathrm{HS}})
    \bigr).
\)
This proves \eqref{mainthm_space'} and \eqref{mainthm_space}.
\end{proof}

The next lemma verifies that the Eulerian profile constructed above
indeed satisfies the generalized Hunter--Saxton equation in the weak
sense.

\begin{lemma}[Verification of the Eulerian weak formulation]
\label{lem_solcheck}
The function
\(
    v(x,t)
    =
    V_{\mathrm{HS}}
    \bigl(\zeta^{-1}(x,t),t\bigr)
\)
is an Eulerian weak solution of \eqref{genHS0} on
\([0,t_*+\veps_{\mathrm{HS}})\). More precisely, for every
\(
    \psi
    \in
    C_c^1
    \bigl(
        \mathbb R\times[0,t_*+\veps_{\mathrm{HS}})
    \bigr),
\)
we have
\begin{equation*}
    \begin{aligned}
    &\int_0^{t_*+\veps_{\mathrm{HS}}}
    \int_{\mathbb R}
    \left[
        v_x\psi_t
        +
        v v_x\psi_x
        +
        \left(1-\frac1\lambda\right)
        v_x^2\psi
    \right]
    \,dx\,dt
    +
    \int_{\mathbb R}
    v_0'(x)\psi(x,0)\,dx
    =
    0.
    \end{aligned}
\end{equation*}
\end{lemma}

\begin{proof}
Set
\(
    S_{\mathrm{HS}}(a,t)
    :=
    v_x(\zeta(a,t),t)
\)
whenever \(D_{\mathrm{HS}}(a,t)\neq0\), and let
\(
    J_{\mathrm{HS}}(a,t)
    :=
    \zeta_a(a,t).
\) 
By \eqref{ux_ricc} and \eqref{eta_a_simple},
\[
    S_{\mathrm{HS}}
    =
    \frac{v_0'}{D_{\mathrm{HS}}},
    \qquad
    J_{\mathrm{HS}}
    =
    D_{\mathrm{HS}}^\lambda
\]
on the regular set.

Introduce
\(
    N_{\mathrm{HS}}
    :=
    S_{\mathrm{HS}}J_{\mathrm{HS}} \) and 
\(
    M_{\mathrm{HS}}
    :=
    S_{\mathrm{HS}}^2J_{\mathrm{HS}}.
\) 
Using the explicit formulas above, we obtain
\[
    N_{\mathrm{HS}}(a,t)
    =
    v_0'(a)
    D_{\mathrm{HS}}(a,t)^{\lambda-1}
\text{ 
and }
    M_{\mathrm{HS}}(a,t)
    =
    v_0'(a)^2
    D_{\mathrm{HS}}(a,t)^{\lambda-2}.
\]
Since \(\lambda\in2\mathbb N\), both expressions extend continuously
through the zero set of \(D_{\mathrm{HS}}\). In particular, they are
locally integrable on
\(
    \mathbb R\times[0,t_*+\veps_{\mathrm{HS}}).
\)

Moreover, since
\[
    \partial_tD_{\mathrm{HS}}(a,t)
    =
    \frac1\lambda v_0'(a),
\]
we have
\begin{equation}\label{HS_N_evolution}
    \partial_tN_{\mathrm{HS}}
    =
    \left(1-\frac1\lambda\right)
    M_{\mathrm{HS}}
\end{equation}
holds across the singular set as well.

Let
\(
    \Psi(a,t)
    :=
    \psi(\zeta(a,t),t).
\)
Since
\(
    \zeta_t(a,t)
    =
    V_{\mathrm{HS}}(a,t)
    =
    v(\zeta(a,t),t),
\)
we have
\[
    \partial_t\Psi(a,t)
    =
    \psi_t(\zeta(a,t),t)
    +
    v(\zeta(a,t),t)
    \psi_x(\zeta(a,t),t).
\]

Since \(a\mapsto\zeta(a,t)\) is strictly increasing, the change of
variables \(x=\zeta(a,t)\) gives
\[
    \begin{aligned}
    &\int_0^{t_*+\veps_{\mathrm{HS}}}
    \int_{\mathbb R}
    \left[
        v_x\psi_t
        +
        vv_x\psi_x
        +
        \left(1-\frac1\lambda\right)v_x^2\psi
    \right]
    \,dx\,dt
    \\
    &\qquad
    =
    \int_0^{t_*+\veps_{\mathrm{HS}}}
    \int_{\mathbb R}
    \left[
        N_{\mathrm{HS}}(a,t)\,
        \partial_t\Psi(a,t)
        +
        \left(1-\frac1\lambda\right)
        M_{\mathrm{HS}}(a,t)\Psi(a,t)
    \right]
    \,da\,dt.
    \end{aligned}
\]
By \eqref{HS_N_evolution}, the integrand on the right-hand side is
\(
    \partial_t
    \bigl(
        N_{\mathrm{HS}}(a,t)\Psi(a,t)
    \bigr).
\)
Hence
\[
    \begin{aligned}
    &\int_0^{t_*+\veps_{\mathrm{HS}}}
    \int_{\mathbb R}
    \left[
        v_x\psi_t
        +
        vv_x\psi_x
        +
        \left(1-\frac1\lambda\right)v_x^2\psi
    \right]
    \,dx\,dt
    =
    -\int_{\mathbb R}
    N_{\mathrm{HS}}(a,0)\Psi(a,0)\,da,
    \end{aligned}
\]
because \(\psi\) is compactly supported in
\(\mathbb R\times[0,t_*+\veps_{\mathrm{HS}})\).

At \(t=0\),
\(
    D_{\mathrm{HS}}(a,0)=1,
    \zeta(a,0)=a,
\) 
and therefore
\(
    N_{\mathrm{HS}}(a,0)=v_0'(a),
    \Psi(a,0)=\psi(a,0).
\) 
Thus
\[
    \begin{aligned}
    &\int_0^{t_*+\veps_{\mathrm{HS}}}
    \int_{\mathbb R}
    \left[
        v_x\psi_t
        +
        vv_x\psi_x
        +
        \left(1-\frac1\lambda\right)v_x^2\psi
    \right]
    \,dx\,dt
    +
    \int_{\mathbb R}
    v_0'(x)\psi(x,0)\,dx
    =
    0.
    \end{aligned}
\]
This establishes the desired weak formulation.
\end{proof} 
Combining the preceding results, we complete the proof of Theorem~\ref{mainthm1}.

\section{The cuspon-anticuspon interaction}\label{sec_HS_global}

In the previous section, we analyzed the solution profile only for
\(t\in(t_\ast,t_\ast+\veps_{\mathrm{HS}})\), during which the solution contains a
single anticuspon--cuspon pair. We now extend our analysis to the global-in-time
regime. This extension is possible because, for the generalized Hunter--Saxton equation,
local minima of \(v_x\) are preserved along the Lagrangian flow. 
For later
times \(t>t_\ast+\veps_{\mathrm{HS}}\), depending on the initial data \(v_0\), new
cuspons may be generated, while existing cuspons may disappear.

Throughout this section, we assume the hypotheses of
Theorem~\ref{thm_collapse}, and in particular
Assumption~\ref{ass:HS-global}. Recall that
\[
\mathcal A:=\{a\in\mathbb R:\ v_0''(a)=0\},
\qquad
\mathcal A^-:=\{a\in\mathcal A:\ v_0'(a)<0\}.
\]
	If \(\mathcal A^-\neq\{a_0\}\), then (G1) and (H3) imply that there exists
	\(a_1\in\mathcal A^-\) such that
	\[
	v_0'(a_1)=\min_{a\in\mathcal A^-\setminus \{a_0\}}v_0'(a).
	\]
	If \(\mathcal A^-=\{a_0\}\), we use the convention \(a_1=\infty\).  
Then, we can find the specific $\veps_{\mathrm{HS}}>0$ that was introduced in Section~\ref{sec_HScusp}--Section~\ref{sec_thmpf}. More precisely, 
\footnote{Note that $\veps_{\mathrm{HS}}\geq 0$, since
	\begin{equation*}
		-\frac{\lambda}{v_0'(a_1)}-t_*\geq -\frac{\lambda}{v_0'(a_0)}-t_*=0.
	\end{equation*}
	Then combining with Assumption~\ref{ass:HS-main}\textup{(H2)}, for $t\in (t_*,t_*+\veps_{\mathrm{HS}})$, we see that there exists exactly two solutions 
	\begin{equation*}
		a_-(t)<a_0<a_+(t), \qquad t\in (t_*,t_*+\veps_{\mathrm{HS}}),
	\end{equation*}
	such that $\zeta_a(a_{\pm}(t),t)=0$.}
\begin{equation*}
	\veps_{\mathrm{HS}}:=	\begin{cases}
		-\dfrac{\lambda}{v_0'(a_1)}-t_*,
		& \text{if } a_1 \in\mathbb{R},\\[1.2em]
		\infty,
		& \text{if } a_1=\infty.
	\end{cases}
\end{equation*}

If \(\mathcal A^-\neq\{a_0\}\), then for times beyond \(t_*+\veps_{\mathrm{HS}}\)
the same cuspon-emergence mechanism as in Lemma~\ref{HS_blowup} may occur
again near the next critical label \(a=a_1\). Thus a new cuspon--anticuspon pair can
emerge in addition to the original one. We now describe how two such pairs evolve relative to each other as time passes. The following
lemma formulates this mechanism for a general initial profile \(v_0\).
\begin{lemma}[Annihilation of an adjacent cuspon-anticuspon pair] \label{lem:HS_cusp_annihilation} Let $a',a''\in \mathcal{A}^-$ be two consecutive local minimum labels of $v_0'$, and suppose that there exists $\overline{a}\in\mathcal{A}$ such that $ a'<\overline{a}<a''$, and \[ v_0^{(3)}(a')>0,\qquad v_0^{(3)}(\overline{a})<0,\qquad v_0^{(3)}(a'')>0. \] In other words, $\overline{a}$ is the unique local maximum label of $v_0'$ between $a'$ and $a''$. 
	For each $a\in\mathcal{A}^-$, let $\beta_{a,+}(t)$ and $\beta_{a,-}(t)$ denote the Lagrangian labels of the cuspon and anticuspon, respectively,  formed on either side of $\zeta(a,t)$ after the corresponding blow-up time $ T(a):=-\lambda/v_0'(a)$.
%
	 \footnote{The labels $a_{\pm}(t)$ introduced in Section~\ref{sec_HScusp} are precisely those corresponding to $\beta_{a_0,\pm}(t)$.}
	
	If $ v_0'(\overline{a})<0$,  set
	\[
	t_{\mathrm{ann}}:=T(\bar a)
	=-\frac{\lambda}{v_0'(\bar a)}.
	\]
	Then, the two inner cuspon labels satisfy \[ \lim_{t\nearrow t_{\mathrm{ann}}}\beta_{a',+}(t) = \lim_{t\nearrow t_{\mathrm{ann}}}\beta_{a'',-}(t) = \overline{a}, \]
	and
	\begin{equation*}
		\lim_{t\nearrow t_{\mathrm{ann}}} v\bigl(\zeta(\beta_{a',+}(t),t),t\bigr) = \lim_{t\nearrow t_{\mathrm{ann}}} v\bigl(\zeta(\beta_{a'',-}(t),t),t\bigr) = v\bigl(\zeta(\overline{a}, t_{\mathrm{ann}}),t_{\mathrm{ann}}\bigr).
	\end{equation*} 
	Moreover, for every \(t>t_{\mathrm{ann}}\), there is no singular label
	in \((a',a'')\). In particular, \(v_x(\cdot,t)\) has no blow-up point in
	\(\zeta((a',a''),t)\).
	
	On the other hand, if \[ v_0'(\overline{a})\geq 0, \] then these two inner cuspons do not annihilate in finite positive time.
\end{lemma}
\begin{proof}
	 By extending the argument of Lemma~\ref{HS_blowup} globally in time, for every $a\in \mathcal{A}^-$, a cuspon pair emerges on both sides of $\zeta(a,t)$ after the corresponding blow-up time  $T(a):=-\lambda/v_0'(a)$. We denote the associated Lagrangian labels by $\beta_{a,\pm}(t)$. Then, for every \(t>\max\{T(a'),T(a'')\}\), both anticuspon--cuspon pairs are
	 present, with associated labels
	 \(\beta_{a',\pm}(t)\) and \(\beta_{a'',\pm}(t)\).
	 
	 Immediately after $t=\max\{T(a'),T(a'')\}$, $\beta_{a',+}(t)<\overline{a}<\beta_{a'',-}(t)$ will hold. We next determine whether these two inner labels can meet at \(\overline{a}\).
	 For this, we first note that
	 since $\beta_{a,\pm}(t)$ is defined to satisfy $D_{\mathrm{HS}}(\beta_{a,\pm}(t),t)=0$, it holds that $v_0'(\beta_{a,\pm}(t))=-\lambda/t$. Differentiating this identity gives 
	 \begin{equation*}
	 	v_0''(\beta_{a,\pm}(t))\beta_{a,\pm}'(t)=\frac{\lambda}{t^2}> 0.
	 \end{equation*}
	 Hence, as long as the corresponding branches remain non-degenerate, we have
	 \begin{equation}\label{apm_glob}
	 	\begin{split}
	 		\beta_{a,+}'(t)> 0 , \quad  \beta_{a,-}'(t)< 0,
	 	\end{split}
	 \end{equation}
	 as $v_0''(\beta_{a,-}(t))<0$ and $v_0''(\beta_{a,+}(t))>0$. Therefore, the distance between the two inner labels $\beta_{a',+}(t)$ and $\beta_{a'',-}(t)$ decreases as time increases. 
	 
	 Suppose that \(v_0'(\bar a)<0\), and set
	 \[
	 t_{\mathrm{ann}}:=T(\bar a)
	 =-\frac{\lambda}{v_0'(\bar a)}.
	 \]
	 For every \(t\in (\max\{T(a'),T(a'')\},t_{\mathrm{ann}})\), 
	 we have
	 \[
	 \max\{v_0'(a'), v_0'(a'')\}<-\frac{\lambda}{t}<v_0'(\bar a).
	 \]
	 Since \(v_0'\) is strictly increasing on \([a',\bar a]\) and
	 strictly decreasing on \([\bar a,a'']\), the equation $v_0'(a)=-\lambda/t$ 
	 has a unique solution \(\beta_{a',+}(t)\) in \((a',\bar a)\), and a unique solution 
	 \(\beta_{a'',-}(t)\) in
	 \((\bar a,a'')\).
	Together with \eqref{apm_glob}, this shows that
	\(\beta_{a',+}(t)\) is increasing and bounded above by \(\bar a\).
	 Let
	 \[
	 \ell_+:=\lim_{t\nearrow t_{\mathrm{ann}}}\beta_{a',+}(t)
	 \in [a',\bar a].
	 \]
	 Passing to the limit in
	 \[
	 v_0'(\beta_{a',+}(t))=-\frac{\lambda}{t},
	 \]
	 we obtain
	 \[
	 v_0'(\ell_+)
	 =
	 -\frac{\lambda}{t_{\mathrm{ann}}}
	 =
	 v_0'(\bar a).
	 \]
	 By the strict monotonicity of \(v_0'\) on \([a',\bar a]\), it follows
	 that \(\ell_+=\bar a\). Hence,
	 \begin{equation}\label{ann_1}
	 	\lim_{t\nearrow t_{\mathrm{ann}}}\beta_{a',+}(t)=\bar a.
	 \end{equation}
	 Similarly,
	 \begin{equation}\label{ann_2}
	 	\lim_{t\nearrow t_{\mathrm{ann}}}\beta_{a'',-}(t)=\bar a.
	 \end{equation}
	 Since \(T\) attains its maximum on \([a',a'']\) at \(\bar a\),
	 the equation \(T(a)=t\) has no solution in \((a',a'')\) for
	 \(t>t_{\mathrm{ann}}\). Hence, there is no singular label in \((a',a'')\) for \(t>t_{\mathrm{ann}}\). In other words, the two inner labels meet at \(t=t_{\rm ann}\) and then disappear.
	 
	 	By \eqref{ann_1}--\eqref{ann_2},
	 \[
	 (\beta_{a',+}(t),t),\,
	 (\beta_{a'',-}(t),t)
	 \longrightarrow
	 (\bar a,t_{\mathrm{ann}})
	 \qquad\text{as } t\nearrow t_{\mathrm{ann}}.
	 \]
	 Since \(V_{\mathrm{HS}}\) is continuous at
	 \((\bar a,t_{\mathrm{ann}})\), we obtain
	 \[
	 \begin{aligned}
	 	\lim_{t\nearrow t_{\mathrm{ann}}}
	 	v(\zeta(\beta_{a',+}(t),t),t)
	 	&=
	 	\lim_{t\nearrow t_{\mathrm{ann}}}
	 	V_{\mathrm{HS}}(\beta_{a',+}(t),t) =V_{\mathrm{HS}}(\bar a,t_{\mathrm{ann}}) =v(\zeta(\bar a,t_{\mathrm{ann}}),t_{\mathrm{ann}}).
	 \end{aligned}
	 \]
	 The same argument applies to \(\beta_{a'',-}(t)\).
	 Therefore, we conclude that 
	 \begin{equation*}
	 	\lim_{t\nearrow t_{\mathrm{ann}}}v(\zeta(\beta_{a',+}(t),t),t)=\lim_{t\nearrow t_{\mathrm{ann}}}v(\zeta(\beta_{a'',-}(t),t),t)=v(\zeta(\overline{a}, t_{\mathrm{ann}}),t_{\mathrm{ann}}).
	 \end{equation*}
	 
	 If \(v_0'(\bar a)\ge 0\), then for every finite \(t>0\),
	 \[
	 -\frac{\lambda}{t}<0\le v_0'(\bar a).
	 \]
	 Hence, for every finite \(t>\max\{T(a'),T(a'')\}\), the two branches
	 remain respectively in \((a',\bar a)\) and \((\bar a,a'')\), and
	 therefore cannot meet in finite positive time.
\end{proof}

\begin{remark}[Change of regularity at cusp annihilation]
	\label{rem:HS_regularity_annihilation}
	Suppose that two adjacent inner cuspon and anticuspon annihilate at the local maximum
	label $\overline{a}$ of $v_0'$, where
	\[
	v_0'(\overline{a})<0,
	\qquad
	v_0''(\overline{a})=0,
	\qquad
	v_0^{(3)}(\overline{a})<0.
	\]
	By Lemma~\ref{lem:HS_cusp_annihilation}, the two inner labels
	meet at time $t_{\mathrm{ann}}
	:=
	-\lambda/ v_0'(\overline{a})$, 
	and the two branches meet at the point 
	$\bigl(\zeta(\overline{a},t_{\mathrm{ann}}),
	v(\zeta(\overline{a},t_{\mathrm{ann}}),t_{\mathrm{ann}})\bigr)$.

	We now examine the regularity of the resulting singular point. By
	\eqref{eta_a_simple}, we have
	\[
	\zeta(a,t_{\mathrm{ann}})
	-
	\zeta(\overline{a},t_{\mathrm{ann}})
	=
	\int_{\overline{a}}^a
	D_{\mathrm{HS}}(\alpha,t_{\mathrm{ann}})^{\lambda}\,d\alpha.
	\]
	Since $\textstyle D_{\mathrm{HS}}(\overline{a},t_{\mathrm{ann}})
	=
	1+\frac{t_{\mathrm{ann}}}{\lambda}v_0'(\overline{a})
	=
	0$ and $\textstyle \partial_\alpha D_{\mathrm{HS}}(\overline{a},t_{\mathrm{ann}})
	=
	\frac{t_{\mathrm{ann}}}{\lambda}
	v_0''(\overline{a})
	=
	0$, 
	Taylor expansion around $\overline{a}$ gives
	\[
	\begin{aligned}
		D_{\mathrm{HS}}(\alpha,t_{\mathrm{ann}})
		&=
		\frac{1}{2}
		\partial_\alpha^2D_{\mathrm{HS}}(\overline{a},t_{\mathrm{ann}})
		(\alpha-\overline{a})^2
		+
		O\bigl(|\alpha-\overline{a}|^3\bigr)
		=
		\frac{t_{\mathrm{ann}}}{2\lambda}
		v_0^{(3)}(\overline{a})
		(\alpha-\overline{a})^2
		+
		O\bigl(|\alpha-\overline{a}|^3\bigr).
	\end{aligned}
	\]
	Consequently, as $a\to\overline{a}$, 
	\[
	\left|
	\zeta(a,t_{\mathrm{ann}})
	-
	\zeta(\overline{a},t_{\mathrm{ann}})
	\right|
	=
	C_{\zeta}|a-\overline{a}|^{2\lambda+1}+o(|a-\overline{a}|^{2\lambda+1}),
	\qquad
	C_{\zeta}>0.
	\]
	Moreover, using $\partial_\alpha V_{\mathrm{HS}}(\alpha,t)
	=
	v_0'(\alpha)D_{\mathrm{HS}}(\alpha,t)^{\lambda-1}$, 
	we obtain, as $a\to\overline{a}$,
	\[
	\left|
	v(\zeta(a,t_{\mathrm{ann}}),t_{\mathrm{ann}})
	-
	v(\zeta(\overline{a},t_{\mathrm{ann}}),t_{\mathrm{ann}})
	\right|
	=
	C_v|a-\overline{a}|^{2\lambda-1}+o(|a-\overline{a}|^{2\lambda-1}),
	\qquad
	C_v>0.
	\]
	Therefore, following the argument of Lemma~\ref{reg_lem}, we have
	\begin{equation*}
		0<\lim_{a\rightarrow\overline{a}}\frac{
			\left|
			v(\zeta(a,t_{\mathrm{ann}}),t_{\mathrm{ann}})
			-
			v(\zeta(\overline{a},t_{\mathrm{ann}}),t_{\mathrm{ann}})
			\right|
		}{
			\left|
			\zeta(a,t_{\mathrm{ann}})
			-
			\zeta(\overline{a},t_{\mathrm{ann}})
			\right|^{
				\frac{2\lambda-1}{2\lambda+1}
			}
		}
		<\infty.
	\end{equation*}
	Hence, two adjacent inner cuspon and anticuspon, each having local regularity
	$C^{\lambda/(\lambda+1)}$, annihilate to form a new singular point at $\zeta(\overline{a},t_{\mathrm{ann}})$, 
	with local regularity $C^{\frac{2\lambda-1}{2\lambda+1}}$. 
	Since $\textstyle \frac{2\lambda-1}{2\lambda+1}
	<
	\frac{\lambda}{\lambda+1}$,
	the annihilation of two inner singularities produces a loss of regularity at
	the collision point.
\end{remark}

\begin{figure}[H]
	\centering
	\includegraphics[width=0.6\linewidth]{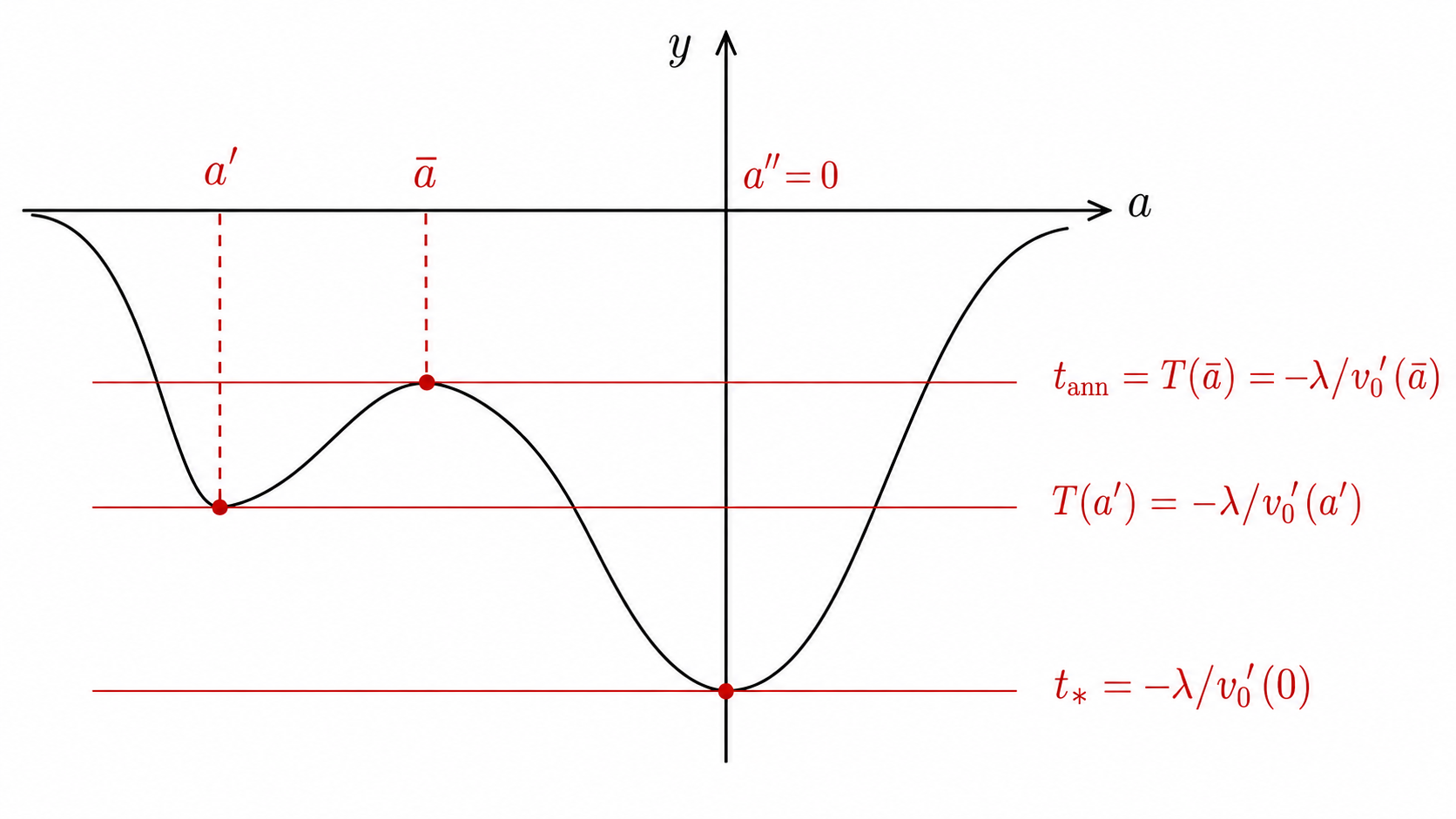}
	\caption{Comparison of the graphs of $v_0'$ and $y=-\lambda/t$.}
	\label{D=0}
\end{figure}

\begin{figure}[H]
	\centering
	\includegraphics[width=0.8\linewidth]{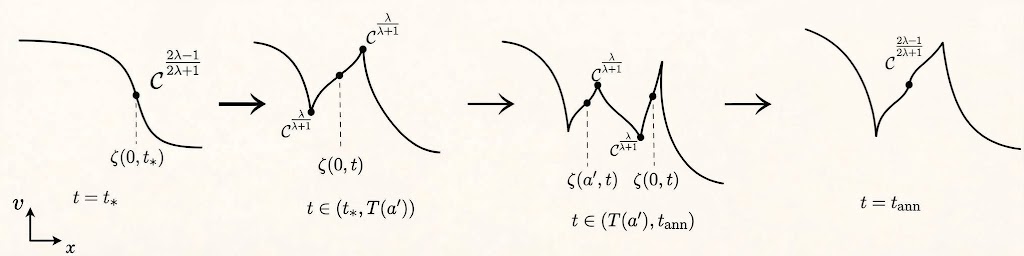}
	\caption{Global-in-time dynamics of $v$ for the initial data in Figure~\ref{D=0}.}
	\label{hs_glob}
\end{figure}

We now examine the annihilation process described above for the initial data given in Figure~\ref{D=0}.
Since $a', \overline{a}, a''$ in Figure~\ref{D=0} satisfy the assumption in Lemma~\ref{lem:HS_cusp_annihilation}, 
cuspon--anticuspon pairs are generated at the times $t_*$ and $T(a')$, while the two adjacent inner singularities annihilate at $t_{\mathrm{ann}}:=T(\overline{a})$. The global dynamics associated with the initial data \(v_0\)
in Figure~\ref{D=0} are given in Figure~\ref{hs_glob}.

\begin{remark}
We expect that this cuspon--anticuspon annihilation phenomenon also occurs for the Camassa--Holm equation.
	For conservative solutions, it has been observed that nonlinear superpositions
	of traveling waves, namely peakons and troughs (antipeakons), of the form
	\[
	u(x,t)=\sum_{j=1}^n p_j(t)\exp\bigl(-2|x-x_j(t)|\bigr),
	\]
	where \(p_j(t)\) and \(-\dot{x}_j(t)\) are positive for peakons and negative
	for antipeakons, pass through one another while preserving their individual
	shapes when they collide \cite{BSS}. Thus, the collision of a cuspon and an
	anticuspon exhibits a type of development that is entirely different from that
	of peakon interactions.
\end{remark}

\section*{Acknowledgements}
B.K. was supported by the National Research Foundation of Korea(NRF) grant funded by the Korea government(MSIT)(00560003).

\end{document}